\documentclass[11pt]{ait} 

\usepackage{cite}
\usepackage{amsmath,amssymb,amsfonts}
\usepackage{algorithmic}
\usepackage{graphicx}
\usepackage{textcomp}
\usepackage{xcolor}

\usepackage{hyperref}

\usepackage{framed}

\year{2021}

\title{%
Смешанная робастность: Анализ систем
с~неопределенными детерминированными
и~случайными параметрами на примере линейных систем
}%
\thanks{Работа выполнена при частичной финансовой поддержке Российского
научного фонда (проект \mbox{№\, 16-11-10015).} 
}

\authors{А.А.~ТРЕМБА, канд.~физ.-мат.~наук (\textrm{atremba@ipu.ru})
\\
(Институт проблем управления им. В.А.~Трапезникова РАН, Москва) }

\abstract{ 
Рассматривается вопрос робастности линейных систем 
с постоянными коэффициентами. 
Существуют методы и инструменты анализа устойчивости
систем со случайными либо детерминированными неопределённостями.
При этом нет подходов для анализа систем, содержащих
оба типа параметрической неопределённости.
В статье проводится классификация типов робастности
и вводится новый тип --- ``смешанная параметрическая робастность'',
включающий несколько вариаций.
Предлагаемые постановки задач о смешанной робастности
могут рассматриваться как промежуточные варианты
между классическим детерминированным и вероятностным 
подходами к робастности.
Перечислены несколько случаев,
в которых поставленные задачи легко решаются.
В общем случае применимы проверки устойчивости 
робастных систем с помощью сценарного подхода, однако эти проверки
могут быть вычислительно сложными.
Для вычисления искомой вероятности устойчивости
предложен простой графо-аналитический подход,
основанный на робастном $D$-разбиении. Этот способ подходит
для случая небольшого числа случайных параметров.
Итоговая оценка вероятности устойчивости вычисляется детерминированным
способом и может быть найдена с произвольной точностью.
Описаны приближенные способы решения поставленных задач.
Приведены примеры и обобщение смешанной робастности на другие типы систем.
}

\begin{document}

\maketitle

{\sl Ключевые слова:}
смешанная робастность, детерминированная робастность, вероятностная робастность, линейные системы, устойчивость.

\section{Введение}
\label{sec:introduction}

В основном в статье рассматриваются динамические системы с постоянными коэффициентами 
и робастность устойчивости относительно неопределённости%
\footnote{%
В общем случае предлагаемые постановки задач
подходят для любых систем, зависящих от параметров,
см. подробнее заключительный раздел~\ref{sec:conclusion}.
Допустимо рассматривать и непараметрическую неопределённость,
но формальное описание соответствующих задач затруднено.
Обобщение смешанной робастности на другие свойства системы или критерии качества, 
помимо устойчивости, тривиально.}.
Во введении намеренно не оговорён точный класс систем,
здесь изложены принципиальные подходы к постановкам и решению задач
робастности.
Далее эти подходы демонстрируются на примере робастной устойчивости.
Без потери общности читатель может рассматривать
экспоненциальную (асимптотическую) устойчивость линейной системы
с одним входом и одним выходом, заданной
дробно-рациональной передаточной функцией.

Основная цель работы состоит в описании нового вида неопределённости ---
смешанной неопределённости и сопутствующей ей робастности 
систем. Этот вид объединяет два типа неопределённости разной природы:
детерминированной и случайной. 

Статья организована следующим образом:
после обзорного введения с описанием различных типов
робастности в разделе~\ref{sec:mixed-robustness} предлагается новый тип
робастности --- смешанная робастность.
Формулировки трёх подтипов задачи о смешанной робастности
приведены в разделе~\ref{sec:problem-statement}
на примере задачи о робастной устойчивости. 
Там же для некоторых из них описан известный способ решения
с помощью сценарного подхода.
В разделе~\ref{sec:main-results} 
для случая маломерного случайного параметра
предлагается использовать графический (графо-аналитический) подход,
основанный на технике робастного $D$-разбиения. 
В подразделе~\ref{sec:approximations} предложены процедуры приближённого решения.
Примеры приведены в разделе~\ref{sec:examples}. 
В заключении приведён ряд направлений развития поставленных задач.

\subsection{Терминология робастности}
Широко используемый в теории управления термин ``робастность''%
\footnote{Англоязычный термин robustness (robust) 
трактуется ещё более разнообразно. На русском языке 
термин робастность по отношению к системам управления 
был систематически использован 
в монографии \cite{polyak-shcherbakov2002}.}
можно понимать по-разному  в зависимости
от контекста и области применения.
Как правило, он используется в качестве характеристики
или свойства (робастные регуляторы, робастные 
системы, робастное управление, робастная стабилизация, робастное оценивание и др.),
но в контексте статьи этот термин будет рассмотрен как самостоятельный,
поскольку интерес представляют разновидности робастности.

Введённый в статистике Боксом и популяризированный Хубером
в контексте ``устойчивости'', ``стабильности'' и ``грубости'' оценок,
термин стал использоваться и в теории управления.

Формирование классической теории робастности состоялась в 80-х гг. XX в. 
В первую очередь исследовалась робастная устойчивость полиномов и матриц,
см. например, дискуссию по проблеме робастности,
состоявшуюся на 11-м Международном конгрессе по автоматическому управлению
IFAC \cite{1992}.
Позднее к ним добавилась $H_\infty$-теория, $M$--$\Delta$ и $\mu$-анализ,
обобщающие различные виды неопределённости.
Систематический обзор результатов представлен в 
\cite{weinmann1991},
\cite{barmish1993},
\cite{bhattacharyya-etal1995}.

В русском языке термин установился не сразу.
Понятие ``грубости'', соответствующее переводу английского ``robust'',
было введено в лексикон теории управления Андроновым в немного другом контексте
(ему соответствует обратное свойство --- ``хрупкость'').
Стало ясно, что нужно слово, описывающее
сохранение свойств при изменении системы или, если смотреть с другой стороны, 
сохранение свойств \emph{множества систем}.

В итоге появился новый термин --- ``робастность'', который 
активно использовал в своих работах Я.З.~Цыпкин и его ученики 
\cite{polyak-tsypkin1991}, \cite{ershov1978}.
Популяризации термина во многом способствовала монография
\cite{polyak-shcherbakov2002}, по сути, закрепившая его.
Вместе с этим для обозначения задач робастности продолжают использоваться 
синонимы ``грубость'', ``нечувствительность'' и др. \cite{braverman-rosonoer1991}.

Каким образом можно выделить класс задач, связанных с робастностью? 
Неформально ``робастностью'' называется свойство нечувствительности
(системы) по отношению к неопределённости.
Далее приведены две устоявшиеся характеристики робастности, а затем ---
два основных типа неопределённости: детерминированная и вероятностная.
На их основе предлагается ввести новый, смешанный тип робастности, 
также удовлетворяющий указанным характеристикам.

Во-первых, ``робастность'' определяется относительно 
неопределённости, которая может быть параметрической 
или непараметрической.
Непараметрическая включает в себя ``немоделируемую динамику'',
описываемую в частотной области ($H_\infty$-теория, QFT), 
структурную неопределённость (например, $M$--$\Delta$ и $\mu$-анализ) и даже внешние помехи.
Часто неопределённость задаётся посредством номинальной системы
и присутствующего в ней возмущения%
\footnote{Здесь термин ``возмущение''
используется не для описания внешнего воздействия,
а как изменение параметров или иных характеристик системы.
Исторически способность системы противодействовать
внешним возмущениям или подавлять внешние шумы 
с \emph{заданными, определёнными} характеристиками 
редко называется робастностью.
Вместе с этим способность системы сохранять заданное свойство 
при \emph{разных} характеристиках внешних сигналов
называться робастностью может (робастность относительно вариативности этих характеристик).
Граница между этими случаями
размыта, устоявшегося общепринятого понимания применимости
термина ``робастность'' нет. 
}.

В статье рассматривается параметрическая статическая неопределённость.
Истинное значение параметра считается неизвестным, 
но зафиксированным в течение времени.
Тем самым можно говорить о \emph{параметрической робастности}.
В этом случае система и её характеристики
явно или неявно зависят от $n$-мерного вещественного вектора 
$q \in \mathbb{R}^n$,
комплекснозначная неопределённость рассматривается аналогично.

Классическое детерминированное описание неопределённости как множества
позже было дополнено вероятностной постановкой,
с помощью случайных величин.
В подразделах~\ref{sec:deterministic-robustness} и \ref{sec:probabilistic-robustness}
эти типы параметрической неопределённости и связанные с ними задачи описаны
подробнее.

Во-вторых, робастность определяется относительно 
\emph{желаемого/требуемого свойства}
системы. Помимо привычной робастной устойчивости, таким свойством может
быть выполнение заданного критерия качества, удовлетворение
фазовым ограничениям или даже сохранение определённых 
статистических характеристик системы \cite{weinmann1991}.
Можно рассуждать о робастности ``такого-то свойства'' по отношению 
к ``такой-то неопределённости''.

Упомянем, что в статистике и анализе данных, где возник термин
``робастность'', он понимается в более узком смысле: 
как чувствительность/нечувствительность к отклонениям распределений 
и выбросам \cite{ershov1978}.

\subsection{Задачи робастности}
Задача \emph{робастного синтеза}
состоит в том, чтобы сделать систему ``робастной'', 
т.е. удовлетворить конкретный критерий для конкретной неопределённости. 
Например, надо найти регулятор заданной структуры, стабилизирующий
систему, параметры которой известны в некоторых пределах,  и т.п.
Полученный регулятор называется ``робастным'', в данном случае
робастно стабилизирующим.

В основе задачи синтеза лежит задача \emph{робастного анализа},
которая состоит в проверке заданного критерия (например, устойчивости)
для заданной неопределённости. 
Далее ``задачей о робастности''
называется именно задача робастного анализа и её вариации. 

Для постановки задачи анализа существенно, как именно проверяется требуемый критерий в зависимости
от типа неопределённости.
Эта проверка может быть бинарной (выполняется/не выполняется),
шкалированной или вероятностной.
Рассмотрим эти задачи подробнее на примере анализа устойчивости.

\subsubsection{``Классический'' робастный параметрический анализ \\
(бинарный, гарантирующий подход)}
\label{sec:deterministic-robustness}

Параметрическую неопределённость 
удобно задавать параллелепипедом (брусом), эллипсоидом,
шаром в какой-либо норме 
или иным заданным множеством  $Q \subset \mathbb{R}^n$. 
Например, характеристический полином системы может
иметь интервальные коэффициенты, которые и являются неопределёнными параметрами.
Данное множество обычно содержит \emph{номинальное} 
значение параметра $q_0 \in Q$, задающего номинальную систему.
Напомним, что в статье рассматривается
статическая параметрическая неопределённость,
т.е. параметр, будучи неизвестным, не меняется с течением времени.

Задача робастной устойчивости состоит в проверке
устойчивости системы \emph{для каждого} элемента из множества 
возможных значений (множества неопределённости):
\begin{equation} \label{eq:deterministic-robustness}
\begin{array}{c}
\text{Параметрическая система} \\
\text{``робастно устойчива''}
\end{array}
 = 
\begin{array}{c}
\text{Система устойчива} \\
\text{для всех } q \in Q
\end{array}
.
\end{equation}

Этот подход реализует
концепцию гарантированной устойчивости, 
принимающий во внимание все ``наихудшие'' варианты (worst-case)
неопределённости.
Будем называть подобный вид робастности 
\emph{детерминированной робастностью},
а задачу --- задачей о детерминированной робастности.
Результат анализа бинарный: либо робастность есть, либо её нет.

\begin{figure}
\begin{center}
\includegraphics[width=6.4cm]{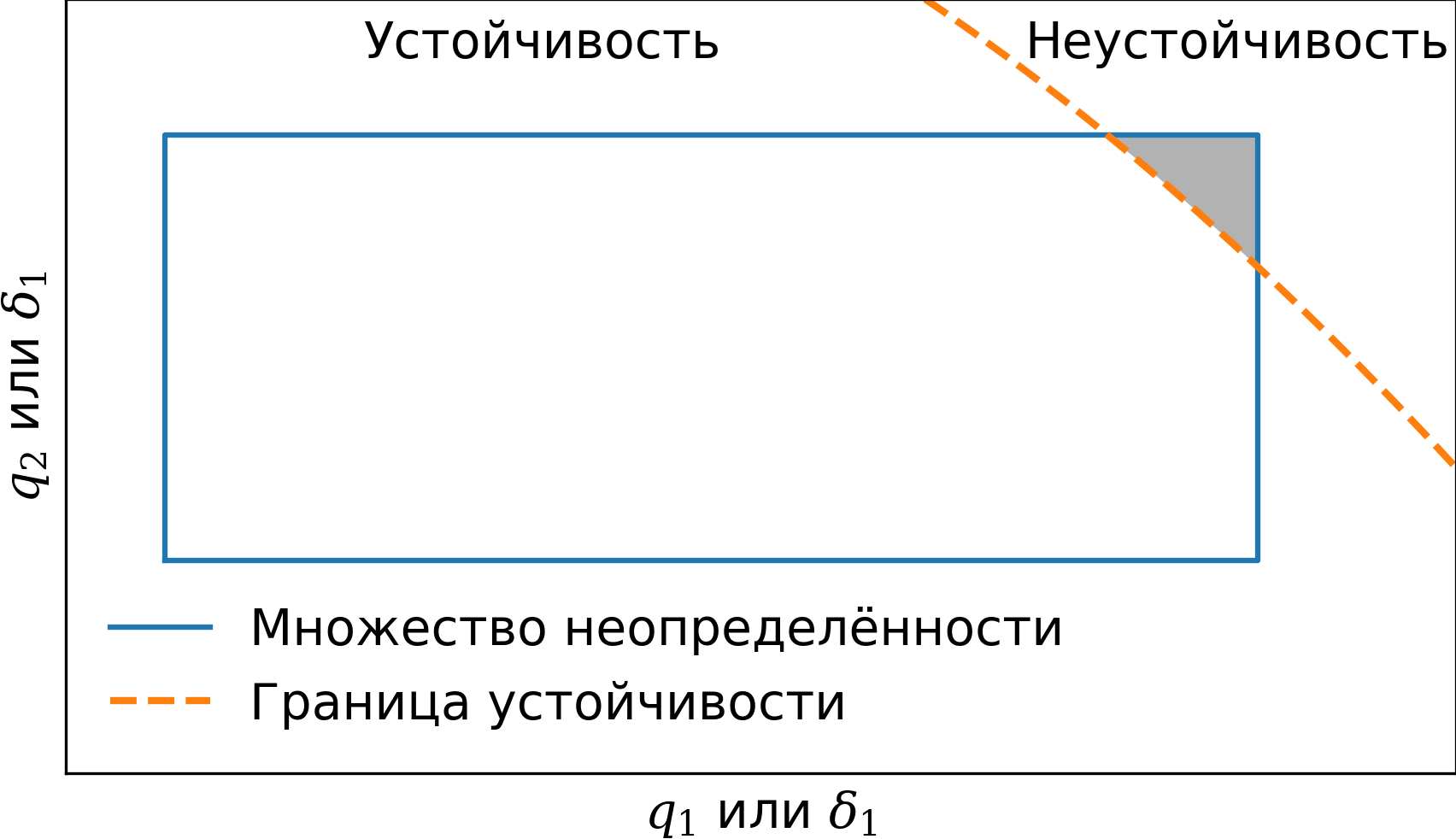}
\caption{Иллюстрация детерминированного и вероятностного подхода.
Параметры из тёмной части прямоугольника неопределённости
соответствуют неустойчивым системам.}
\label{fig:robustness}
\end{center}
\end{figure}

В качестве примера можно рассмотреть прямоугольник на рис.~\ref{fig:robustness} 
как множество неопределённых параметров $Q$ и устойчивость в качестве критерия. 
Каждое значение параметра --- точка в двумерном пространстве $(q_1, q_2)$ ---
соответствует либо устойчивой, либо неустойчивой системе.
Пусть множества ``устойчивых'' и ``неустойчивых'' параметров
разделены некоторой ``пограничной'' линией, отсекающей
часть прямоугольника-неопределённости.
Поскольку часть множества $Q$ принадлежит области неустойчивых параметров,
то решением задачи детерминированной робастной устойчивости
будет ``Нет'': система не робастно устойчива.

Проверка детерминированной робастности (а тем более
решение задачи робастного синтеза)
в большинстве случае является принципиально сложной задачей.
Для анализа устойчивости существует ряд успешно разрешённых
вопросов, таких как Теорема о вершинах, Рёберная теорема или знаменитый результат
Харитонова о сведении анализа устойчивости
интервального семейства полиномов к анализу всего четырёх полиномов
\cite{kharitonov1978},\cite{bhattacharyya-etal1995}. 
Однако эти случаи являются скорее исключением, чем правилом.
Для большинства интервальных случаев эта задача является NP-полной,
например для анализа робастности по отношению
к неопределённости-многограннику может потребоваться 
перебор всех вершин этого множества,
а число вариантов перебора экспоненциально растёт с размерностью 
\cite{nemirovskii1993},\cite{poljak-rohn1993}.


Помимо сложности решения задачи о классической робастности,
критике подвергалась сама постановка: почему
множество допустимых неопределённых параметров выбирается 
тем или иным, каков размер этого множества? Например, если 
неопределённость задана шаром, то как выбирать его радиус?
Отметим, что в классическом робастном анализе есть
важная задача о радиусе робастной устойчивости (robustness margin), 
состоящая в определении максимального размера множества неопределённости
заданной формы, сохраняющего устойчивость
\cite{polyak-shcherbakov2002},\cite{bhattacharyya-etal1995}.

Для полноты картины следует упомянуть подход, 
описывающий робастность не-бинарного критерия качества,
например перерегулирования.
В зависимости от размера множества неопределённости $Q$
значение критерия будет ``деградировать'' (ухудшаться при увеличении
размера) \cite{tempo-etal2013}. Кривая
деградации в зависимости от размера неопределённости 
даёт полезную информацию для принятия
решений. Такой подход назван шкалированной, модулированной робастностью
(modulated robustness) \cite{garatti-campi2013}. 
Для бинарного же критерия, такого как устойчивость, переход от устойчивости
к неустойчивости ``мгновенный'' и определяет радиус робастной устойчивости.

Кроме того, множества неопределённости часто формируются
исходя из грубых верхних и нижних оценок внутренних параметров,
что неизбежно влечёт консерватизм решения.
Ещё одной проблемой классического робастного анализа является
встречающаяся ``сосульчатая'' (icicle) геометрия множества неустойчивых параметров 
\cite{garatti-campi2013}. 
В этом случае есть параметры, близкие к номинальным, делающие систему неустойчивой,
и радиус робастной устойчивости близок к нулю.
Но при этом объём множества таких ``неустойчивых'' параметров
может быть пренебрежимо мал
по сравнению со всем множеством неопределённости. 

Таким образом, будучи несомненно полезным и важным подходом к исследованию 
и синтезу систем управления, классический, гарантирующий подход 
к робастности имеет свои ограничения.

\subsubsection{Вероятностный робастный параметрический анализ \\
(случайные параметры)} \label{sec:probabilistic-robustness}

Для преодоления указанных проблем был предложен вероятностный подход. 
В этом подходе неопределённость в системе считается \emph{вероятностной},
представленной случайным вектором $\delta \in \mathbb{R}^m$ \cite{tempo-etal1996}.
Этот вектор имеет некоторое распределение на множестве $\Delta \subset \mathbb{R}^m$,
являющемся носителем (функции) распределения.
Для простоты%
\footnote{Строго говоря, распределение задаётся тройкой 
$(\Delta, \mathcal{F}, \mathcal{P})$ 
с множеством элементарных событий $\Delta$, 
сигма-алгеброй его подмножеств $\mathcal{F}$ 
и вероятностной мерой $\mathcal{P} : \mathcal{F} \rightarrow [0, 1]$,
но для краткости избыточные обозначения опущены. 
Далее мера обозначается символом $\mu$.
Используется традиционная мера Лебега, а
обозначение распределения ассоциируется с множеством $\Delta$.
Рассматриваемые множества неопределённости предполагаются измеримыми.
Можно рассматривать как абсолютно непрерывные распределения, 
заданные функцией плотности распределения, так и дискретные распределения
либо их смесь.
В контексте робастного анализа с вероятностной точки зрения важен их носитель
и (кумулятивная) функция распределения.} 
обозначим такой вид неопределённости вместе с его распределением
как $\delta \sim \Delta$.
Существенно, что в данном подходе \emph{разрешается} потеря устойчивости
(как правило, небольшая), которая описывается в терминах теории вероятности.

Задача \emph{вероятностного анализа} робастной устойчивости
состоит в нахождении вероятности устойчивости $p^*$
по отношению к случайной неопределённости $\delta$:
\begin{equation} \label{eq:probabilistic-problem}
p^* = \mathrm{Prob}_{\delta \sim \Delta} \big[\text{Система устойчива для данного } \delta \big].
\end{equation}

Эта вероятность равна мере множества ``устойчивых'' параметров,
т.е. значений, при которых система устойчива.
Будем называть этот вид устойчивости по отношению к ``случайной'' 
неопределённости \emph{вероятностной робастностью}.
Это уже не бинарное свойство, оно сопровождается 
своим значением вероятности.

Вероятность логического утверждения, зависящего от случайной переменной $\delta$,
понимается как мера множества, на котором это утверждение истинно. Она равна
интегралу Лебега от индикаторной функции этого утверждения, принимающей
значения 0, если утверждение ложно, и 1, если истинно.
Будем предполагать, что утверждение алгоритмически разрешимо (индикаторная
функция вычислима в точках носителя и на его подмножествах).
Автор намеренно не использовал термин ``стохастическая робастность'',
так как этот термин более широкий и используется чаще 
по отношению к устойчивости систем со случайными процессами.
Следует особо отметить, что система не обладает стохастической динамикой,
но часть её параметров является реализацией конечномерной 
случайной величины $\delta$.

Задача о вероятностной робастности \eqref{eq:probabilistic-problem} 
связана с детерминированной робастностью \eqref{eq:deterministic-robustness}
несколькими путями.

Во-первых, следуя подходу \cite{barmish-lagoa1997}, предлагается
\emph{методически заменять} детерминированную неопределённость $q \in Q$
случайной неопределённостью  $\delta$ той же размерности,
\emph{равномерно распределенной} на множестве $Q$ (тем самым полагая $\Delta = Q$).
Таким образом, ``трудная'' задача о детерминированной робастности 
(включая задачи синтеза) заменяется полиномиально сложной \cite{bai-etal1997}.
Это связано с тем, что можно эффективно
вычислять вероятность устойчивости \eqref{eq:probabilistic-problem}
методом Монте-Карло (см. подраздел~\ref{sec:scenario} с подробностями). 
Очевидно, что новая задача имеет другое содержание,
но её решение может быть интерпретировано в контексте исходной 
детерминированной задачи.
В случае единичной вероятности устойчивости $p^* = 1$ 
результат можно трактовать%
\footnote{С точностью до множеств нулевой вероятности. В зависимости от задачи
эта разница может быть ничтожной либо, наоборот, --- существенной.}
как детерминированную робастность по отношению к множеству
$Q$, совпадающему с $\Delta$.
Если же вероятность устойчивости меньше единицы, то это соответствует
отсутствию детерминированной робастной устойчивости. Полученный результат
может привести к разумной корректировке постановки исходной задачи, например,
в виде рекомендации уменьшить множество неопределённых параметров.

Во-вторых, неопределённость может быть вероятностной изначально,
т.е. описывать случайность значений параметров. 
В этом случае распределение уже не обязательно равномерное.
Согласно классическому подходу к робастности, в качестве множества неопределённости
выбирается либо весь носитель случайных параметров, либо его доверительное
множество высокой вероятности (множество квантили). 
В обоих случаях в постановку вносится 
консерватизм, сопровождаемый вопросами
выбора уровня квантили и выбора соответствующего доверительного множества.
В то же время исходная вероятностная задача анализа 
может быть эффективно решена напрямую, и ответом является вероятность устойчивости. 
Более того, для успешного решения достаточно иметь генератор выборки 
распределения случайных параметров/систем, а не его аналитическое описание \cite{tempo-etal2013}.
Этот же подход используется для вероятностных задач синтеза, см. \cite{calafiore-campi2006}.

Ключевая идея --- использовать или вводить случайность неопределённых параметров --- 
получила должную популярность: 
основная монография \cite{tempo-etal2013} была переиздана спустя 8 лет.
Подход получил широкое распространение и в приложениях,
поскольку он строго обосновывает идеи \emph{рандомизации} 
при анализе и синтезе систем управления.

Рассмотрим тот же рис.~\ref{fig:robustness}.
Пусть прямоугольник является носителем $\Delta$ 
с равномерным распределением случайного параметра $\delta$.
Тогда решение задачи вероятностного робастного анализа \eqref{eq:probabilistic-problem}
даётся в точности отношением площади незаштрихованной части прямоугольника
к его общей площади. Наряду с вероятностью устойчивости 
удобно использовать \emph{вероятность потери устойчивости} $1 - p^*$.
Для равномерного распределения она равна доле заштрихованной части 
прямоугольника, т.е. 2~\%.

Оценка вероятности устойчивости \eqref{eq:probabilistic-problem} 
может быть выполнена несколькими способами.

Во-первых, если дано описание множества устойчивых параметров,
то меру этого множества можно найти явно. Для равномерного
распределения это задача геометрической вероятности:
надо найти отношение объёмов устойчивых параметров в носителе
к объёму самого носителя $\Delta$. Эта же величина может быть вычислена
с помощью объёма множества неустойчивых параметров 
(см. формулу \eqref{eq:probability-by-ratio} далее).
Сами объёмы могут быть найдены аналитически или численно.

Во-вторых, оценка вероятности $p^*$ может быть сделана
эмпирически с помощью случайной выборки 
на множестве неопределённости согласно
сопутствующему ему распределению.
Для равномерного распределения отношение 
подмножеств оценивается отношением числа точек, попавших в каждое
из подмножеств. 
Этот подход является самым известным и мощным инструментом
для решения вероятностных задач анализа и синтеза.
В теории управления он известен как 
\emph{сценарный подход} и использует идею \emph{рандомизации},
см. подробнее подраздел~\ref{sec:scenario}.
В его основе лежит метод типа Монте-Карло, 
включающий вычисление необходимого размера выборки для достижения
заданного \emph{качества} оцениваемой вероятности 
\cite{tempo-etal1996},\cite{calafiore-campi2006}.
Обширная коллекция рандомизированных алгоритмов для решения
вероятностных задач управления приведена в монографии
\cite{tempo-etal2013} и в её библиографическом списке.
Сценарный подход может быть использован и для решения задач синтеза,
которые в статье не рассматриваются.

Отметим, что существуют способы формирования выборки для 
сложных и неявно заданных множеств неопределённости.
В них используются реализации метода Монте-Карло в марковских цепях, 
организованные как случайное блуждание ``биллиардного'' типа (Billiard Walk)
или ``удар-отход'' (Hit-and-Run) \cite{gryazina-polyak2014}. 

Есть и обратная зависимость между вероятностной и детерминированной
робастностью. Задачу о вероятностной робастности 
можно приближённо решать детерминированными методами.
Для фиксированного значения вероятности $p$ выбирается множество квантили
$Q_p : \mathrm{Prob}_{\delta \sim \Delta} [\delta \in Q_p] = p$,
и для него проверяется робастная устойчивость. Если система 
оказывается устойчивой (в детерминированном смысле относительно $Q_p$), 
то вероятность робастной устойчивости не меньше $p$.
Подбор уровня и множества квантили 
остаётся во многом эвристической процедурой, во многом из-за
``сосульчатого'' множества неопределённых параметров.

\section{Смешанная робастность}
\label{sec:mixed-robustness}

В данном разделе неформально вводится концепция
\emph{смешанной робастности}, совмещающей два вышеупомянутых типа робастности.
Строгие формулировки приведены в разделе~\ref{sec:problem-statement}.

Рассмотрим линейную систему с постоянными коэффициентами, 
имеющую и детерминированную параметрическую неопределённость $q$,
и случайную $\delta$. Параметры 
принадлежат своим множествам $q \in Q \subset \mathbb{R}^n$, 
$\delta \in \Delta \subseteq \mathbb{R}^m$.
Случайный вектор $\delta$ распределён на множестве $\Delta$ и обозначается как
$\delta \sim \Delta$. 
Параметры распределения здесь несущественны и опущены.

Начнём с простого случая: пусть распределение случайного параметра $\delta$
не зависит от значения параметра $q$ (также неопределённого).
Анализ устойчивости таких систем включает в себя
и вероятностную, и гарантирующую части.
Соответственно задача является смесью детерминированной \eqref{eq:deterministic-robustness} 
и вероятностной \eqref{eq:probabilistic-problem} задач.
Она состоит в нахождении вероятности
\begin{equation} \label{eq:informal-mixed-robustness}
\mathrm{Prob}_{\delta \sim \Delta}[\text{Система устойчива для данного } \delta
\text{ и всех } q \in Q].
\end{equation}

Формально поставленная задача может называться
задачей о ``вероятностной робастности''.
В конце концов, \emph{это и есть} задача вероятностного анализа
устойчивости типа \eqref{eq:probabilistic-problem},
сформулированная с более сложным утверждением
``Система \emph{робастно} устойчива (по $q$ на множестве $Q$) для данного $\delta$''
вместо ``Система устойчива для данного $\delta$''.
Её разрешимость с помощью сценарного подхода
и связь с известными способами решения обсуждается 
в подразделе~\ref{sec:scenario}.
Более того, эта задача подпадает под обобщённую формулировку
задачи~6.4 из \cite{tempo-etal2013}
о вычислении (вероятностного) показателя качества системы.
Вместе с этим существуют более удобные методы вычисления
вероятности \eqref{eq:informal-mixed-robustness},
использующие особенность постановки.

Во избежания двусмысленности термина ``робастность''
предпочтительно использовать термин 
\emph{смешанная параметрическая робастность}.
Этот термин подчёркивает не только параметрическую форму неопределённости, 
но и существенную роль разных
типов неопределённости: случайной и не случайной.
Поскольку далее рассматриваются только параметрические неопределённости,
уточнение ``параметрическая/ие'' опускается.

В общем случае параметры этих двух видов могут быть зависимыми 
или подчинёнными друг другу, см. раздел~\ref{sec:problem-statement} статьи. 
Соответствующие постановки задач являются также новыми.
В разделе \ref{sec:main-results} предложен ряд подходов к решению,
специфичных для задач смешанной робастности.

Отметим, что довольно часто встречаются задачи,
в которых детерминированная параметрическая
робастность смешивается с вероятностной \emph{динамической} неопределённостью
в шуме/возмущениях/ограничениях, см.
\cite{fan-etal1991},
\cite{jayasuriya1993},
\cite{reinelt-ljung2001},
\cite{vitus-etal2016},
\cite{silvadeaguiar-etal2018}.
Такие неопределённости также известны под именем структурированно-неструктурированных
возмущений (structured-unstructured perturbation)
\cite{chapellat-etal1990}.
Другим примером смешивания с зависимыми случайными
и не случайными переменными являются задачи оптимизации с квантильными ограничениями,
но в них преследуется иная цель.
В то же время комбинация случайных и не случайных неопределённых 
параметров изучалась редко или неявно, 
для других задач или в другом контексте,
в частности для упрощения
задачи синтеза регулятора \cite{fujisaki-etal2008}
или для оценки некоторого распределения \cite{wu-etal2016}.

Далее понятие смешанной робастности проводится на примере анализа
линейных систем управления с постоянными коэффициентами,
содержащих смешанную (гибридную, совместную) параметрическую
неопределённость.

\subsection{Мотивирующие факторы для введения
смешанной параметрической робастности}

Деление между детерминированными и случайными параметрами
возникает в нескольких случаях.

Во-первых, в системе эти параметры могут иметь различную природу.
Например, случайный параметр $\delta$ 
представляет неопределённость в регуляторе
(неточность реализации), 
в то время как параметр $q$ отвечает за неопределённость 
в управляемой системе либо наоборот.
В \cite{fujisaki-etal2008} подобное разделение
было сделано искусственно для не робастной задачи синтеза
регулятора. 
В этом случае часть параметров регулятора случайно выбирались
по методу Монте Карло, а области устойчивости по остальным параметрам
легко вычислялись. Основной идеей являлся случайный, 
рандомизированный подбор регуляторов,
этот подход не был связан с анализом робастности.

Во-вторых, выделение нескольких параметров из общего числа варьируемых параметров
 полезно для более глубокого исследования системы.
Этот выбор может быть сделан вручную для оценки так называемой ``модулированной'' 
(шкалированной) робастности по отношению к выбранным параметрам
либо для оценки потери качества при изменении множества неопределённости
\cite{garatti-campi2013},\cite{tempo-etal2013}.
Случайность выделенных параметров может быть навязана
искусственно \cite{bai-etal1997}.
С одной стороны, эта случайность моделирует нежёсткие ограничения,
а с другой --- отражает априорные представления о неопределённости.
В зависимости от задачи случайный параметр можно рассматривать, например, 
как равномерно или нормально распределенный, 
с несимметричным распределением, представленный смесью распределений и т.д.
Такой подход помогает не только формализовать априорные
представления о неопределённости, но также помогает 
убедиться в отсутствии \emph{хрупкости},
т.е. негативной чувствительности системы к малым изменениям выбранных параметров.
А параметры, найденные с помощью оптимизации (без учёта робастности),
часто именно такие, подробнее об этом эффекте см. \cite{keel-bhattacharyya1997}.

В-третьих, предложенная постановка может возникнуть, 
если некоторая часть неопределённых параметров
ограничена ``жёстко'' ($q$), 
а другие параметры ограничены ``мягко'' ($\delta$). 
Действительно, даже если использовать вероятностный подход,
не обязательно ожидать, что \emph{все} неопределённые параметры 
представимы случайными величинами.
Этот случай включает в себя нечёткие неопределённости,
которые, например, могут быть смоделированы 
как сумма интервального и случайного значений.
При этом привлекать специальный аппарат
нечёткого исчисления не требуется.

В-четвёртых, задача о смешанной робастности является
результатом \emph{частичной} рандомизации, при которой
только некоторые из исходно детерминированных
неопределённых параметров заменяются случайными.

Основным результатом статьи является постановка задач
для смешанной параметрической робастности, 
включающей различные типы неопределённости.
В случае нескольких случайных параметров 
для решения этой задачи предлагается использовать
графо-аналитический подход.

Предложенное понятие смешанной робастности по сложности заполняет 
промежуток между детерминированными задачами и
изложенными в \cite[гл. 6]{tempo-etal2013} 
идеями вероятностного анализа этих задач.
Тем самым можно найти компромисс между сложностью
(вложенной) задачи о детерминированной робастности
и ``неуверенностью'' (т.е. случайностью исхода)
решения вероятностной задачи о робастности.

Подчеркнём схожесть и отличия между предложенными
задачами о смешанной робастности и вероятностным
вычислением показателя качества из \cite[гл. 6]{tempo-etal2013}.
В то время как задача~6.1 в \cite{tempo-etal2013} поставлена
относительно произвольного критерия, 
в приведённых примерах предполагается, что
случайно распределенный параметр $\delta$ --- это \emph{единственная}
неопределённость.
Затем эта неопределённость исследуется вероятностными методами
с учётом детерминированных ограничений на неё.
Таким образом, одна из задач ставится как 
поиск вероятности%
\footnote{Значение показателя качества в постановке заменено на
критерий устойчивости для краткости.
Упомянутая глава содержит много других идей анализа, например
изучение деградации критерия качества (при случайных параметрах) и т.п.,
см. также похожий подход к ``шкалированной'' робастности в \cite{garatti-campi2013}.
}
$$
\mathrm{Prob}_{\delta \sim \Delta}\big[P(s, \delta) \text{ гурвицев \textbf{и} } \delta \in Q\big].
$$

В предложенной же задаче о смешанной робастности \eqref{eq:informal-mixed-robustness} 
критерий, параметризованный случайной переменной $\delta$,
включает внутреннюю робастность по отношению
к \emph{другой неопределённости} $q$. 
Тем самым простейшая задача о смешанной робастности
одновременно охватывает две \emph{различные} группы ``робастности'':
внутреннюю по отношению к параметру $q$ 
и внешнюю по отношению к параметру $\delta$.
Если бы $\delta$ был тоже детерминированным
(в контексте классической робастности), либо параметр $q$ был бы случайным,
то в итоге задача свелась бы к классической детерминированной
постановке \eqref{eq:deterministic-robustness} или 
к полностью вероятностной постановке \eqref{eq:probabilistic-problem}
соответственно.

Помимо рассмотренного примера
независимых случайных и неслучайных параметров,
возникают два новых типа смешанной робастности:
случайные параметры могут зависеть от неслучайных 
и наоборот --- множество неопределённости детерминированных 
параметров само может быть случайным, подробнее см. раздел~\ref{sec:problem-statement}.

Смешивание разных типов параметрической робастности и их взаимосвязь 
в задачах автоматического управления являются новыми.
С другой стороны, идея вычисления вероятности попадания
случайного параметра $\delta$ в ``хорошее'' или ``плохое''
множество хорошо известна, см., например, \cite[раздел 6.2]{tempo-etal2013}. 
Трудность состоит в удобном описании этих множеств.
Для случая маломерного случайного параметра эта задача решена 
в разделе~\ref{sec:main-results} с помощью робастного $D$-разбиения.

\section{Постановка задач о смешанной робастности} 
\label{sec:problem-statement}

Рассмотрим замкнутую линейную систему с постоянными коэффициентами
без запаздывания%
\footnote{Расширения на другие типы систем представлены 
в разделе~\ref{sec:conclusion}.} 
с двумя типами параметров:
детерминированными $q \in Q \subset \mathbb{R}^n$ и
случайными $\delta \sim \Delta \subseteq \mathbb{R}^m$.
Случайная величина $\delta$ распределена на множестве $\Delta$
с некоторой вероятностной мерой (равномерной, нормальной и т.д.).
Ещё раз подчеркнём, что предложенные постановки имеют смысл не только
для задач устойчивости, но и для любых желаемых свойств системы.

Устойчивость рассматриваемой системы определяется корнями
характеристического полинома, коэффициенты которого зависят
от параметров
\[
P(s, q, \delta) = a_0(q, \delta) s^k + ... + a_k(q, \delta).
\]
Для устойчивости непрерывной системы этот полином должен быть
гурвицевым, т.е. все его корни $s_i$ должны иметь
отрицательную вещественную часть $\mathrm{Re}\, s_i < 0$.
Для устойчивости дискретной системы аналогичный полином
$P(z, q, \delta)$ должен быть шуровским, т.е. все корни
должны находиться внутри единичного круга $|z_i| < 1$.
Без потери общности далее рассматривается только 
непрерывный случай, и устойчивость системы
отождествляется со свойством Гурвица её характеристического
полинома.

Задача анализа \emph{смешанной параметрической робастности}
состоит в том, чтобы найти вероятность устойчивости
по случайной переменной $\delta$ 
при наличии неопределённости $q$.
Поскольку в общем случае случайные и детерминированные параметры
могут быть зависимыми друг от друга,
приведём три постановки задачи о смешанной робастности.

В самом простом варианте случайная и детерминированная неопределённости 
полностью независимы, т.е. задача имеет форму \eqref{eq:informal-mixed-robustness}.
Пусть это и самый простой способ смешивания, но он является новым. 
Здесь неопределённые параметры равноправны.

\begin{framed}
I. {\sl{З\,а\,д\,а\,ч\,а\,}} о независимой смешанной параметрической робастности 
($Q$--$\Delta$~робастность).
{\itshape

Пусть заданы случайный параметр $\delta \sim \Delta \subseteq \mathbb{R}^m$
и множество неопределённости $Q \subset \mathbb{R}^n$, а также
характеристический полином системы с постоянными коэффициентами
$P(s, q, \delta)$.
Требуется найти вероятность устойчивости
\begin{equation} \label{eq:mixed-robustness-independent}
\!\!\!p^* = \mathrm{Prob}_{\delta \sim \Delta}
\Big[ P(s, q, \delta) \text{ гурвицев для всех } q \in Q \Big].
\end{equation}
}
\end{framed}

В статье эта и последующие значения вероятности называются 
\emph{вероятностями смешанной робастности} или 
\emph{вероятностями смешанной робастной устойчивости}.

Далее, множество возможных параметров $Q$
для детерминированной неопределённости
может зависеть от случайного параметра $\delta$, т.е.,
строго говоря, само быть случайным.
Тем самым детерминированная неопределённость подчинена случайной.

\begin{framed}
II. Смешанная {\sl{з\,а\,д\,а\,ч\,а\,}} робастности с зависимостью первого типа
($Q(\Delta)$ робастность).
{\itshape

Пусть заданы случайный параметр $\delta \sim \Delta \subseteq \mathbb{R}^m$,
значения которого связаны с множествами
неопределённости $Q(\delta) \subset \mathbb{R}^n$,
и характеристический полином системы с постоянными
коэффициентами $P(s, q, \delta)$. 
Требуется найти вероятность устойчивости
\begin{equation} \label{eq:mixed-robustness-Q(Delta)}
\!p^* \!\!=\! \mathrm{Prob}_{\delta \sim \Delta}
\!\Big[\! P(s, \delta, q) \text{ гурвицев для всех } q \in Q(\delta) \Big].
\end{equation}
}
\end{framed}

Отметим, что характеристический полином зависит не только от случайного
множества неопределённости $q$, но и от самого случайного параметра
$\delta$. 
Можно рассматривать форму записи задачи только со случайным множеством
(без неопределённости $\delta$ в характеристическом полиноме)
$\Big[\! P(s, \bar{q}) \text{ гурвицев для всех } \bar{q} \in \bar{Q}(\delta) \Big]$.
Легко показать, что эти постановки эквивалентны.

Наконец, случайный параметр $\delta$ (по существу --- его распределение 
и носитель) может зависеть от значения детерминированной неопределённости.
Эту постановку можно рассматривать как поиск ``наихудшего''
варианта неопределённости. При этом неопределённость 
понимается в расширенном смысле и включает неопределённость в распределении
$\Delta(q)$.

\begin{framed}
III. Смешанная {\sl{з\,а\,д\,а\,ч\,а\,}} робастности с зависимостью второго типа
($\Delta(Q)$ робастность).
{\itshape

Пусть заданы множество неопределённых параметров 
$Q \subset \mathbb{R}^n$
и семейство случайных параметров $\delta \sim \Delta(q) \subseteq \mathbb{R}^m$,
параметризованное с помощью $q \in Q$.
Для характеристического полинома системы с постоянными коэффициентами
$P(s, q, \delta)$ 
найти \emph{гарантированную} вероятность устойчивости
\begin{equation} \label{eq:mixed-robustness-Delta(Q)}
p^* = \inf_{q \in Q} \mathrm{Prob}_{\delta \sim \Delta(q)}
\Big[\! P(s, q, \delta) \text{ гурвицев} \Big].
\end{equation}
}
\end{framed}

Эта постановка охватывает параметрические семейства распределений
как частный случай, при котором меняются только параметры, но не носитель.
Для удобства положим, что указанный минимум достигается (например, если
$Q$  компактно, а распределение задано функцией распределения, равномерно
непрерывной по параметру $q$).
Вместе с этим предложенная постановка более общая, так как устойчивость системы 
может зависеть не только от \emph{параметризованного распределения} случайного параметра
($\min_{q \in Q} \mathrm{Prob}_{\delta \sim \Delta(q)} [P(s, \delta) \text{ гурвицев}]$),
но детерминированный параметр может также влиять на устойчивость непосредственно, 
например, частью компонент, к примеру
$\min_{(q^a, q^b) \in Q} \mathrm{Prob}_{\delta \sim \Delta(q^a)} [P(s, q^b, \delta) \text{ гурвицев}]$.

Как частный случай, в третий тип попадает задача с фиксированным распределением,
при котором неопределённость входит только в сам характеристический полином:
\[
\min_{q \in Q} \mathrm{Prob}_{\delta \sim \Delta} [P(s, q, \delta)  \text{ гурвицев}].
\]
Несмотря на то что в такой постановке оба типа параметров независимы, 
эта задача принципиально отличается от 
задачи о $Q$--$\Delta$ робастности \eqref{eq:mixed-robustness-independent}.
Её также можно интерпретировать как задачу
о нахождении \emph{гарантированного ``среднего'' показателя устойчивости},
понимая вероятность устойчивости
$\mathrm{Prob}_{\delta \sim \Delta} [P(s, q, \delta)  \text{ гурвицев}]$ 
как ``усреднённое значение
устойчивости'' по множеству $\Delta$. 
Схожий подход к минимизации не бинарного критерия описан 
в \cite{garatti-campi2013}.
Отметим, что попутно может быть определено ``наихудшее''
с точки зрения устойчивости значение параметра $q^* \in \mathrm{Arg} \min_{q \in Q} (\cdots)$.
Такое значение параметра может быть не единственным, даже если минимум достижим.
Если минимум не достигается, возможна постановка задачи о нахождении предельной точки
или локализации области ``наихудших'' параметров.

С точки зрения разрешимости поставленных задач самой простой, по-видимому,
является задача о $Q$--$\Delta$ робастности 
с независимыми неопределённостями \eqref{eq:mixed-robustness-independent}.
Методы, подходящие для решения задачи 
о $Q(\Delta)$ робастности \eqref{eq:mixed-robustness-Q(Delta)}, 
также подходят и для неё.
Задача о $\Delta(Q)$ робастности с обратной зависимостью имеет принципиально другую 
природу. По-видимому, она наиболее трудоёмка (за исключением дискретного случая, см.
подраздел~\ref{sec:discrete}).
Для неё в подразделе~\ref{sec:approximations} предлагается использовать
приближенный метод.

Далее перечислены несколько способов решения
поставленных задач о смешанной робастности с помощью известных методов.
В подразделе~\ref{sec:approximations} предложены
методы приближённого решения.
В разделе~\ref{sec:main-results} задачу анализа о независимой $Q$--$\Delta$ 
робастности предлагается решать в два этапа с помощью 
вспомогательной задачи синтеза.

\subsection{Аналитическое решение}
Для первых двух задач о смешанной робастности ($Q$--$\Delta$ и $Q(\Delta)$) 
можно формально описать внутренний критерий 
с помощью индикаторной функции устойчивости
\begin{equation} \label{eq:indicator-function-robustness}
F(\delta) = I[P(s, q, \delta) \text{ робастно устойчив}].
\end{equation}
Здесь робастность подразумевается либо для фиксированного множества
неопределённости $Q$, либо для параметризованного $Q(\delta)$.
Для разных значения $\delta \in \Delta$ 
вычисление этой функции сводится к различным
задачам о детерминированной робастности.
В этом случае решение задачи о смешанной робастности
\eqref{eq:mixed-robustness-independent} 
записывается как
\begin{equation} \label{eq:exact-solution}
p^* = \mathop{\int}_{\delta} F(\delta) d \mu,
\end{equation}
т.е. как мера множества \emph{робастно стабилизирующих} параметров $\Delta_{good} \ni \delta$.
Здесь использована мера Лебега $\mu$, лежащая в основе распределения
случайного параметра.

Однако функция $F$ редко доступна в явном виде,
поскольку её вычисление включает детерминированную задачу 
о робастной устойчивости \eqref{eq:deterministic-robustness},
которая может быть трудоёмка.

Аналогично характеризуется прямой способ решения задачи $\Delta(Q)$-робастности.
Если вероятность устойчивости в параметризованных вероятностных задачах робастности
$p^*(q) = \mathrm{Prob}_{\delta \sim \Delta(q)}
\big[\! P(s, q, \delta) \text{ гурвицев} \big]$ имеет явный вид, 
то решение можно получить с помощью задачи оптимизации $p^* = \min_{q \in Q} p^*(q)$,
в общем случае невыпуклой и негладкой.

\subsection{Сценарный подход {\rm (}выборка методом Монте-Карло{\rm)}}
\label{sec:scenario}
Следуя идее рандомизации в вероятностном анализе,
для оценки вероятности \eqref{eq:exact-solution}
может быть использован сценарный подход (выборочный, эмпирический)
\cite{tempo-etal1996}.
Надо выбрать $N$ случайных образцов $\delta^i, i = 1,..., N$, 
согласно распределению параметра $\delta$ 
и получить долю устойчивых систем, которую можно записать
с помощью индикаторной функции устойчивости как
\[
\widehat{p}^* = \frac{1}{N} \sum_{i=1}^{N} F(\delta^i).
\]
В этом случае неравенство
\begin{equation} \label{eq:chernoff-bounds-on-p*}
|p^* - \widehat{p}^*| < \varepsilon
\end{equation}
будет выполнено с вероятностью
$1 - 2e^{-2 \varepsilon^2 N} \equiv 1 - \theta$
для любой \emph{точности} (accuracy) $\epsilon \in (0, 1)$.
Эта вероятность характеризуется уровнем достоверности (confidence level) $\theta$. 
Достаточного числа образцов, необходимого для достижения
(высокой) вероятности $1 - \theta$
с данной (также высокой) точностью $\varepsilon$ 
оценки \eqref{eq:exact-solution}
даётся верхней оценкой Чернова \cite{tempo-etal1996}:
\begin{equation} \label{eq:chernoff-formulae}
N = \Big\lceil \frac{1}{2 \varepsilon} \ln \frac{2}{\theta}\Big\rceil.
\end{equation}
Например, для достижения точности $\varepsilon = 0{,}01$
и уровня достоверности $\theta = 10^{-7}$
нужны $84057$ образцов.

Формирование выборки является отдельной задачей, легко решаемой
для произвольных непрерывных маломерных распределений (например, методом обращения или отсечения), 
некоторых распределений векторов и матриц, 
а также для дискретных распределений \cite[главы 14, 16, 18]{tempo-etal2013}.

Размер выборки может быть значительно уменьшен, если использовать 
свойство выпуклости критерия. Априори ожидать выпуклости индикаторной
функции устойчивости \eqref{eq:indicator-function-robustness} нельзя
из-за подзадачи о робастной устойчивости.
Таким образом, для достоверного вычисления вероятности устойчивости
типичный размер выборки должен быть велик, что особенно затрудняет
применение сценарного подхода для задач смешанной робастности.

Как уже упоминалось в разделе~\ref{sec:introduction}, 
в общем случае каждая подзадача сложна.
Для прямого же применения сценарного подхода
данные задачи должны легко решаться для каждого элемента выборки.
И действительно, такие случаи есть.

Робастная устойчивость полинома с независимыми
интервальными коэффициентами может быть проверена
аналитически с помощью теоремы Харитонова \cite{kharitonov1978}.
В этом случае для образца $\delta^i$ каждый из коэффициентов $a_j(q, \delta^i)$ 
должен \emph{независимо} меняться в интервалах
$\underline{a}_j(\delta^i) \leq a_j(q, \delta^i) \leq \bar{a}_j(\delta^i)$.
Тогда достаточно проверить гурвицевость четырёх ``угловых'' полиномов.
Отметим, что эта проверка верна только для характеристического 
полинома системы непрерывного времени, для дискретного случая
аналогичного результата нет \cite{jury1990}.

Отметим также, что для любого числа образцов полученная оценка
\eqref{eq:chernoff-bounds-on-p*} эмпирическая
и \emph{всегда} будет выполнена не гарантированно,
а лишь с большой долей уверенности, заданной параметром $\theta$.

\subsection{Приближенные методы} \label{sec:approximations}

Задачи о смешанной робастности подразумевают ответ в виде численного значения, 
поэтому альтернативой точному решению могут служить верхние и нижние границы
вероятности устойчивости.

Можно использовать простые
необходимые или достаточные условия робастной устойчивости
для вычисления функции $F(\delta^i)$, 
что приведёт к нижним или верхним границам
вероятности смешанной робастности.
Так, если вместо индикаторной функции устойчивости $F(\delta^i)$
в формуле \eqref{eq:exact-solution} 
использовать индикаторную функцию \emph{необходимого} критерия устойчивости,
то это даст \emph{верхнюю} границу вероятности устойчивости.
Соответственно использование \emph{достаточного} критерия
устойчивости даст \emph{нижнюю} границу вероятности \eqref{eq:exact-solution}.
В контексте сценарного подхода вычисление критериев
может быть намного проще проверки робастной устойчивости, 
при этом оценка вероятности \eqref{eq:chernoff-bounds-on-p*}
будет односторонней,
а необходимые или достаточное условия устойчивости могут быть слишком
консервативными и в результате давать грубые оценки вероятностей
смешанной робастности.

Задачу о $Q(\Delta)$ робастной устойчивости с зависимыми параметрами
можно приближённо решать как задачу с независимыми параметрами.
Во-первых, если итоговое множество ограничено, то 
можно получить верхнюю и нижнюю границу, объединяя или 
пересекая множества неопределённости:
$$
p^*_{II} \geq p^*_{I,\cup} =  \mathrm{Prob}_{\delta \sim \Delta}
\Big[ P(s, q, \delta) \text{ гурвицев для всех } 
q \in Q = \bigcup_{\delta \in \Delta} Q(\delta) \Big],
$$
$$
p^*_{II} \leq p^*_{I,\cap} =  \mathrm{Prob}_{\delta \sim \Delta}
\Big[ P(s, q, \delta) \text{ гурвицев для всех } 
q \in Q = \bigcap_{\delta \in \Delta} Q(\delta) \Big].
$$
Во-вторых, в контексте сценарного подхода для решения задачи 
о $Q(\Delta)$ робастности в качестве нижней оценки 
вероятности устойчивости можно использовать решение
задачи о $Q$--$\Delta$ робастности относительно объединения множеств 
$Q = \bigcup_{i=1,...,N} Q(\delta^i)$.

Для всех трёх типов задач о смешанной робастности применима
идея аппроксимации носителя случайного распределения с помощью
множества квантили, упомянутая в конце подраздела~\ref{sec:probabilistic-robustness}.
Так, если для фиксированного $p$ и множества 
$Q_p : \mathrm{Prob}_{\delta \sim \Delta} [\delta \in Q_p] \geq p$
полином $P(s, q, \delta)$ робастно устойчив в классическом смысле
по отношению к $q \in Q, \delta \in Q_p$, то вероятность робастной
устойчивости \eqref{eq:mixed-robustness-independent} с независимыми
параметрами не меньше $p$.
Для задачи \eqref{eq:mixed-robustness-Q(Delta)} о $Q(\Delta)$ робастности
надо проверять робастную устойчивость характеристического полинома 
по отношению к $q \in \cup_{\delta \in Q_p} Q(\delta), \delta \in Q_p$.
Наконец, этот приём можно использовать для задачи \eqref{eq:mixed-robustness-Delta(Q)} 
о $\Delta(Q)$ робастности, выбрав 
$Q_p : \forall q \in Q, \mathrm{Prob}_{\delta \sim \Delta(q)} [\delta \in Q_p] \geq p$.

\subsection{Робастное $D$-разбиение}

Рассмотрим метод синтеза маломерных систем управления,
который будет полезен для решения задач о смешанной робастности:
задачу о робастной $D$-устойчивости \cite{petrov-polyak1991-rus}. 

Пусть интересующая нас система с постоянными коэффициентами
зависит от параметров $d \in \mathbb{R}^m$
и неопределённости $q \in Q \subset \mathbb{R}^n$. 
Её характеристический полином $P(s, q, d)$.
Предположим, что множество локализации корней $D$ 
в комплексной плоскости также является частью постановки задачи.
Для задачи робастной стабилизации непрерывных систем
это левая комплексная полуплоскость
$D = \{s : \mathrm{Re}\, s < 0\} \subset \mathbb{C}$ и т.д.

Цель задачи синтеза состоит в том, чтобы предоставить такое значение
параметра $d$, при котором все корни характеристического полинома 
$P(s, q, d)$ лежат в множестве $D$ \emph{для всех допустимых неопределённостей}:
\begin{equation} \label{eq:robust-parameteric-problem}
d^* \colon \text{Все корни полинома } P(s, q, d^*) \text{ лежат в } D \text{ для всех } q \in Q.
\end{equation}

Данная задача иногда называется задачей о робастном расположении
корней или задачей синтеза о робастной кластеризации (локализации).
Надлежащим выбором множества $D$ можно решать различные задачи синтеза: \\
\phantom{\hspace{0.5cm}} а) о робастной стабилизации систем непрерывного времени 
($D$ --- левая комплексная полуплоскость), \\
\phantom{\hspace{0.5cm}} б) о робастной стабилизации дискретных систем 
($D$ --- внутренность единичного круга), \\
\phantom{\hspace{0.5cm}} в) о робастной стабилизации с ограниченной колебательностью 
($D$ является сектором в левой комплексной полуплоскости) \\
\phantom{\hspace{0.5cm}} и т.д. 

Есть разные способы решения этой задачи от алгебраических
и оптимизационных до основанных на линейных матричных неравенствах.
В практически каждом из этих подходов выдаётся единственное значение
$d^*$, которое является ``наилучшим'' по отношению 
к дополнительно выбранному критерию.

Альтернативным подходом является описание 
\emph{всего} множества робастно стабилизирующих параметров, 
которые удовлетворяют задаче \eqref{eq:robust-parameteric-problem}:
\[
D_{stab} \!
= \! \{d^* \!:\! \text{ Все корни полинома } P(s, q, d^*) \text{ лежат в } 
D \text{ для всех } q \in Q \}.
\]
Робастное $D$-разбиение является графо-аналитическим
методом и чаще всего применяется для одно- и двухмерных параметров.
Он позволяет неявно изобразить множество устойчивых параметров 
$D_{stab} \subseteq \mathbb{R}^m$,
а точнее, его \emph{границу} \cite{petrov-polyak1991-rus}. 
Эта граница строится с помощью принципа исключения нуля 
и ``заметания'' (перебора по сетке) обобщённых частот 
либо с помощью аналитического вывода уравнений,
см. \cite{tremba2006-rus}, 
\cite{siljak1992}, 
\cite{hwang-etal2010}, 
\cite{mihailescu-stoica-etal2017}
и цитируемые там источники.

Основным ограничением практического применения
робастного $D$-разбиения является низкая размерность 
пространства параметров. Обычно оно может быть сделано 
для $m = 1, 2$ или максимум для трёх параметров.

Снова отметим, что указанная задача является задачей \emph{синтеза}.
Она же может быть ключевым компонентом для решения
задачи \emph{анализа} смешанной робастности 
с независимыми параметрами \eqref{eq:mixed-robustness-independent}.

\section{Основные результаты} \label{sec:main-results}
Для решения задачи о смешанной робастности с независимыми параметрами \eqref{eq:mixed-robustness-independent}
или зависимостью первого типа \eqref{eq:mixed-robustness-Q(Delta)}
предлагается отказ от сценарного подхода 
в части формирования подзадач робастного анализа.
В общем случае для решения задач о смешанной робастности 
с зависимыми параметрами предлагается использовать аппроксимации 
из подраздела~\ref{sec:approximations}).

Наличие дискретной неопределённости (случайной и/или непрерывной)
существенно упрощает решение задач о смешанной робастности.

\subsection{Смешанная робастность с независимыми группами детерминированных и~случайных параметров
{\rm(}$Q$--$\Delta$ робастность{\rm)}}

Идея состоит в том, чтобы рассматривать случайный параметр $\delta$
как детерминированный и для него найти множество ``хороших'', ``робастно устойчивых''
параметров
\[
\Delta_{stab} = \{\delta : F(\delta) = 1\} 
= \{\delta : \text{ система робастно устойчива для всех } q \in Q\}.
\]
Напомним, что искомая вероятность смешанной робастной устойчивости
\eqref{eq:mixed-robustness-independent}
является вероятностной мерой этого множества. Единственным ограничением
подхода является требование его измеримости.

Для нахождения этого множества предлагается использовать\footnote{Этот способ 
особенно хорошо подходит для анализа робастности любого свойства системы, 
определяемого расположением корней характеристического полинома.} 
робастное $D$-разбиение
\emph{в пространстве случайной неопределённости} $\delta \in \mathbb{R}^m$
по отношению к ``жёстко'' ограниченной неопределённости $q \in Q$.
Каждая точка множества $\Delta_{stab} \equiv D_{stab}$ обеспечивает
робастную устойчивость характеристического полинома по отношению
к детерминированной неопределённости $q$.
Затем решение задачи \eqref{eq:mixed-robustness-independent}
вычисляется как (вероятностная) мера этого множества.

Предлагается следующий алгоритм решения задачи о смешанной робастности
\eqref{eq:mixed-robustness-independent}.
Он записан в терминах устойчивости линейной системы с постоянными параметрами,
но тривиально модифицируется на иные свойства систем.

\begin{framed}
\begin{algorithm}
\mbox{\textbf{Решение задачи о смешанной робастности ($Q$--$\Delta$ робастности)}}
\label{alg:independent-Q-Delta}

\begin{enumerate}
\item[{\sl Ш\,а\,г 1.}]
Найти множество робастно стабилизирующих параметров
$\Delta_{stab}$ по отношению к детерминированной
неопределённости $q$
\[
\Delta_{stab} = \{\delta : P(s, q, \delta) \text{ гурвицев для всех } q \in Q\}.
\]
\item[{\sl Ш\,а\,г 2.}]
Для данной случайной неопределённости $\delta \sim \Delta$
вычислить вероятность устойчивости как
\begin{equation} \label{eq:probability-by-good-set}
p^* = \mathrm{Prob}_{\delta \sim \Delta} [\delta \in \Delta_{stab}] = \mu(\Delta_{stab}).
\end{equation}
\end{enumerate}
\vspace{-4mm}
\end{algorithm}
\end{framed}

Подчеркнём, что каждый шаг алгоритма может быть выполнен независимо с помощью
любого подходящего метода.
Например, оценка вероятности \eqref{eq:probability-by-good-set} 
на шаге~2 может быть сделана методом Монте-Карло (по сути, тем же сценарным подходом).
Этот подход отличается от
прямого сценарного подхода из подраздела~\ref{sec:scenario}. 
На шаге~2 каждый случайно выбранный образец $\delta^i$
должен быть проверен на принадлежность конструктивно
описанному множеству $\Delta_{stab}$, 
в то время как в прямом сценарном подходе 
для каждого образца должна решаться подзадача о детерминированной робастности.

Если распределение случайной неопределённости $\delta$ равномерно на $\Delta$,
то геометрическая вероятность \eqref{eq:probability-by-good-set}
вычисляется как отношение объёмов (площадей) множеств
\begin{equation} \label{eq:probability-by-ratio}
p^* 
= \frac{\mathrm{Vol}(\Delta \cap \Delta_{stab})}{\mathrm{Vol}(\Delta)} 
= 1 - \frac{\mathrm{Vol}(\Delta \setminus \Delta_{stab})}{\mathrm{Vol}(\Delta)}.
\end{equation}
При использовании робастного $D$-разбиения граница множества $\Delta_{stab}$ 
задаётся поточечно (численно либо сеткой при аналитическом задании границы), 
тем самым это множество аппроксимируется многоугольником
или набором многоугольников, площадь которых может быть найдена точно.
Уменьшая сетку, эта аппроксимация позволяет
найти площадь $\Delta_{stab}$ с требуемой точностью. 

Для не равномерного распределения задача \eqref{eq:probability-by-good-set} 
о численном нахождении меры множества $\Delta_{stab}$ может быть нетривиальной
даже при аналитическом описании множества.
В одномерном случае задача решается явно (см. подраздел~\ref{sec:scalar-random-uncertainty}), 
в двухмерном и трёхмерном случаях можно использовать перебор одного или двух 
параметров.
Вопрос о точности полученных сеточных аппроксимаций и используемых квадратурных формул 
требует отдельного рассмотрения.

Основным недостатком предлагаемого подхода является низкая размерность
случайного параметра $\delta$. 
Выполнимость первого шага ограничена применимостью техники
анализа пространства параметров, такого как робастное $D$-разбиение,
поскольку множество $\Delta_{stab}$ 
может быть конструктивно построено только в 
одно-, дву- и (максимум) трёхмерных случаях.

Другой особенностью является численная сложность построения 
самого множества $\Delta_{stab}$. Если параметры $q, \delta$ 
входят в коэффициенты характеристического полинома аффинно, 
то робастное $D$-разбиение можно построить численно или аналитически
\cite{petrov-polyak1991-rus},
\cite{tremba2006-rus},
\cite{hwang-etal2010}.
В противном случае построение робастного
$D$-разбиения может быть сложной задачей само по себе.
В таких случаях предлагается использовать внутренние или внешние
приближения множества $\Delta_{stab}$, получая для них соответственно 
нижние и верхние оценки вероятности \eqref{eq:probability-by-good-set}.

Для задачи о смешанной $Q(\Delta)$ робастности (с зависимостью первого типа)
можно применять этот же алгоритм~\ref{alg:independent-Q-Delta}, 
если на первом шаге удаётся найти множество 
робастно стабилизирующих параметров $\Delta_{stab}$ с учётом зависимости:
\[
\Delta_{stab} = \{\delta : P(s, q, \delta) \text{ гурвицев для всех } q \in Q(\delta)\}.
\]

Далее перечислены несколько случаев, в которых задача
о смешанной робастности или её подзадачи могут быть решены явно.

\subsection{Скалярная случайная неопределённость} 
\label{sec:scalar-random-uncertainty}
Пусть случайный параметр системы одномерен: $\delta \in \mathbb{R}^1$.
Если множество $\Delta_{stab}$ непусто, то оно может быть
представлено как объединение непересекающихся отрезков%
\footnote{Первый и последний отрезки могут иметь бесконечную длину, например 
быть $(-\infty, b_1]$.
Также эти отрезки могут быть открытыми или полуоткрытыми,
например для каких-то $j$ быть $[a_j, b_j]$, для других $(a_j, b_j)$ и т.д.
Разница между этими случаями имеет меру нуль
и для абсолютно непрерывного распределения $\delta$ 
не влияет на вычисление вероятности.
Однако в случае дискретного или смешанного распределения эти нюансы необходимо учитывать.}
$(a_i, b_i], \; i = 1, ...$, 
где $a_i < b_i$ и $b_i < a_{i+1}$.

Эти интервалы могут быть найдены разными способами,
например с помощью одномерного робастного $D$-разбиения
\cite{petrov-polyak1991-rus} или
так называемым методом ``разведчика устойчивости'' (Stability Feeler),
предложенным в \cite{matsuda-mori2009}.

Таким образом, система является робастно устойчивой для некоторого
значения $\delta$ тогда и только тогда, когда
существует интервал $i : a_i < \delta \leq b_i$. Значение
индикаторной функции устойчивости \eqref{eq:indicator-function-robustness} 
тривиально вычисляется
$$
F(\delta) = \left\{
\begin{array}{cl}
1, & \text{ если } \exists i : \delta \in (a_i, b_i], \\
0\phantom{,} & \text{ в противном случае.}
\end{array}
\right.
$$

Для найденных отрезков робастной устойчивости
задача о смешанной робастности
\eqref{eq:mixed-robustness-independent} имеет очевидное 
решение \eqref{eq:exact-solution}:
\begin{equation} \label{eq:exact-solution-1d}
p^* = \sum_{i= 1, ...} \big(\mathrm{\bf CDF}(b_i) - \mathrm{\bf CDF}(a_i)\big),
\end{equation}
т.е. является суммой вероятностей попадания
неопределённого параметра $\delta$ в каждый из отрезков $(a_i, b_i]$.
Решение записано с помощью одномерной функции распределения
неопределённого параметра $\delta$
\[
\mathrm{\bf CDF}(a) = \mathrm{Prob}_{\delta \sim \Delta} [\delta \leq a].
\]

\subsection{Дискретная неопределённость} \label{sec:discrete}

Пренебрегая переходными процессами, анализ робастности системы с переключениями 
можно свести к анализу системы с дискретной неопределённостью.
Аналогичным образом конечный набор однотипных систем
может быть записан в виде одной системы с добавлением в неё 
неопределённого дискретного параметра.

Анализ смешанной робастности значительно упрощается 
при наличии дискретной неопределённости любого типа. 
Приёмы, подобные перечисленным далее, можно использовать в случаях,
когда дискретны только некоторые компоненты случайного и/или детерминированного параметров.

\subsubsection{Дискретная детерминированная неопределённость}
Если число различных значений детерминированных параметров конечно и невелико
(например, не превышает $K$), то подзадачи робастного анализа 
$Q$--$\Delta$ и $Q(\Delta)$ робастной устойчивости 
\eqref{eq:mixed-robustness-independent}, \eqref{eq:mixed-robustness-Q(Delta)}
сводятся к проверке устойчивости небольшого числа полиномов.
В этом случае разумно применять прямой сценарный подход из подраздела~\ref{sec:scenario},
решая эту подзадачу для конечного же числа образцов $\delta^i, 
i = 1,...$:
``$P(s, \delta^i, q^{i,k}) \text{ гурвицев для всех } k \in K_i \subseteq {1, ..., K}$''.

Задача $\Delta(Q)$ о смешанной робастности сводится к перебору 
и решению конечного числа подзадач о вероятностной робастности
``$\mathrm{Prob}_{\delta \sim \Delta_k}
\Big[\! P(s, q^k, \delta) \text{ гурвицев} \Big], k = 1, ..., K$''.

\subsubsection{Дискретная случайная неопределённость}

Пусть распределение задано набором значений $\delta^j$
вместе с вероятностями $p_j = \mathrm{Prob}_{\delta \sim \Delta} [\delta = \delta^j], \, j = 1, ...$
Если множество допустимых значений $\delta^j$ 
конечно и имеет небольшой размер, $j = 1, ..., M$, 
то связанная с ним задача о смешанной $Q(\Delta)$ (или $Q$--$\Delta$)
робастности может быть решена явно.

\begin{framed}
\begin{algorithm}
\label{alg:discrete-random}
\textbf{Решение задачи о смешанной $Q(\Delta)$ робастности c дискретной случайной
неопределённостью \eqref{eq:mixed-robustness-Q(Delta)}}

Пусть задана случайная величина $\delta$,
принимающая значения $\delta^j$ с вероятностями $p_j$, 
и характеристической полином системы $P(s, q, \delta)$
зависит также от неопределённого параметра $q \in Q(\delta)$.

Для каждого $j = 1, ..., M$ решить задачу о детерминированной робастности 
и вычислить
\[
\!p^* \!\!=\! \sum_{j=1}^M p_j \cdot I\Big[\! P(s, q, \delta^j) \text{ гурвицев для всех } q \in Q(\delta^j) \Big].
\]
\end{algorithm}
\end{framed}

Технически, этот подход похож на сценарный, но
число ``образцов'' ограничено, а в результате будет получено точное решение.
Анализ распадается на $M$ независимых подзадач о детерминированной робастности
и его сложность зависит от сложности решения этих подзадач.
Как и в случае сценарного подхода, эти подзадачи можно решать параллельно.

Отметим, что этот результат не зависит от размерности
случайной неопределённости $\delta$.

Алгоритм~\ref{alg:discrete-random} применим не только к задаче \eqref{eq:mixed-robustness-Q(Delta)} 
о смешанной робастности c зависимостью, 
но и к задаче \eqref{eq:mixed-robustness-independent}
с независимыми параметрами.

Можно аппроксимировать решение, отбросив ряд ``маловероятных'' 
значений параметра $\delta$, и тем самым получить нижнюю оценку вероятности устойчивости.
Так же можно поступить, если множество дискретных значений бесконечно.

\subsubsection{Дискретная случайная и дискретная детерминированная неопределённость}
Данный случай примечателен тем, что в случае конечного числа значений параметров
все три варианта задачи о смешанной робастности: $Q$--$\Delta$, $Q(\Delta)$ и $\Delta(Q)$
легко решаются.

В задаче о $Q$--$\Delta$ робастности случайный параметр задан набором значений $\delta^j$
вместе с вероятностями $p_j = \mathrm{Prob}_{\delta \sim \Delta} [\delta = \delta^j], \, j = 1, ..., M$,
а детерминированный --- набором $Q = \{q^k, k = 1,..., K\}$.
Решение даётся формулой
\[
\!p^* \!\!=\! \sum_{j=1}^M p_j \cdot \prod_{k=1}^K I\Big[\! P(s, q^k, \delta^j) \text{ гурвицев}\Big].
\]

В задаче о $Q(\Delta)$ робастности детерминированный параметр задан набором с двумя индексами
$\{q^{j, k}, k = 1,..., K_j\}$ для каждого $j = 1, ..., M$. 
Решение даётся аналогичной формулой
\[
\!p^* \!\!=\! \sum_{j=1}^M p_j \cdot \prod_{k=1}^{K_j} I\Big[\! P(s, q^{j, k}, \delta^j) \text{ гурвицев}\Big].
\]

Наконец, в задаче о $\Delta(Q)$ робастности детерминированная неопределённость
задана набором $\{q^k, k = 1,..., K\}$, а случайная -- наборами
значений $\{\delta^{k, j}, j = 1,...,M_k\}, k = 1,..., K$, c вероятностями
$p_{k,j} = \mathrm{Prob}_{\delta \sim \Delta_k} [\delta = \delta^{k, j}], \, j = 1, ..., M_k$. 
В этом случае решение даётся формулой
\[
p^* = \min_{k = 1,...,K} \sum_{j=1}^{M_j} p_{k,j} \cdot I\Big[\! P(s, q^k, \delta^{k,j}) \text{ гурвицев}\Big].
\]

\section{Примеры} \label{sec:examples}

В первых трёх примерах рассматриваются задачи о независимой
смешанной робастности ($Q$--$\Delta$ робастности), 
далее в примерах решаются задачи о зависимой смешанной робастности.
В четвёртом примере рассмотрена задача о зависимости первого типа ($\Delta(Q)$ 
робастность), а в пятом -- зависимость второго типа ($\Delta(Q)$ робастность).
В последних примерах активно используется возможность описать множество устойчивости
анализируемой системы на расширенном множестве неопределённостей $Q \times \Delta$.

\subsection{Скалярная случайная неопределённость} \label{sec:example-1-scalar-random}
Рассмотрим немного изменённый пример~2 из \cite{matsuda-mori2009}
об устойчивости дискретных систем.
Объект управления с передаточной функцией $G(z, K)$,
включающей коэффициент усиления, 
стабилизируется управлением с обратной связью с помощью регулятора $C(z, q_1, q_2)$:
\[
G(z, K) = \frac{K}{1 + 3 z},\;\;\;
C(z, q_1, q_2) = \frac{8{,}5 + 3{,}5 z}{1 + (8 + q_1) z + (6 + q_2) z^2}.
\]
Каждый из неопределённых параметров $q_1, q_2$ находится в интервале $[-1, 1]$.

Задача о смешанной робастности состоит в том, 
чтобы определить вероятность устойчивости,
рассматривая коэффициент усиления обратной связи
$K$ как нормально распределенную случайную величину со средним
$2$ и стандартным отклонением $0{,}1$. 
Устойчивость определяется критерием Шура--Кона (а не Гурвица) 
для характеристического полинома
$$
P(z, q, K) = K(8{,}5 + 3{,}5 z) + (1 + 3 z)\big(1 + (8 + q_1) z + (6 + q_2) z^2\big).
$$

В \cite{matsuda-mori2009} определена робастная устойчивость этой системы
при $K \in [-0{,}11,\, 2{,}21]$, что соответствует первому шагу 
алгоритма~\ref{alg:independent-Q-Delta}.
Формула \eqref{eq:exact-solution-1d} для нормального распределения принимает вид
$$
p^* = \frac{1}{0{,}1 \sqrt{2 \pi}} \mathop{\int}_{-0{,}11}^{2{,}21} e^{-\frac{(K - 2)^2}{2\cdot (0{,}1)^2}} dK 
= \frac{1}{2}\Big( 
\mathrm{erf} \Big(\frac{2{,}21 - 2}{0{,}1 \sqrt{2}}\Big)
- \mathrm{erf} \Big(\frac{-0{,}11 - 2}{0{,}1 \sqrt{2}}\Big)
\Big) = 0{,}982.
$$

В результате получена точная вероятность смешанной робастной устойчивости.

\subsection{Двумерная случайная неопределённость}  \label{sec:example-3-2d-random}

\begin{figure}
\centerline{\includegraphics[width=9cm]{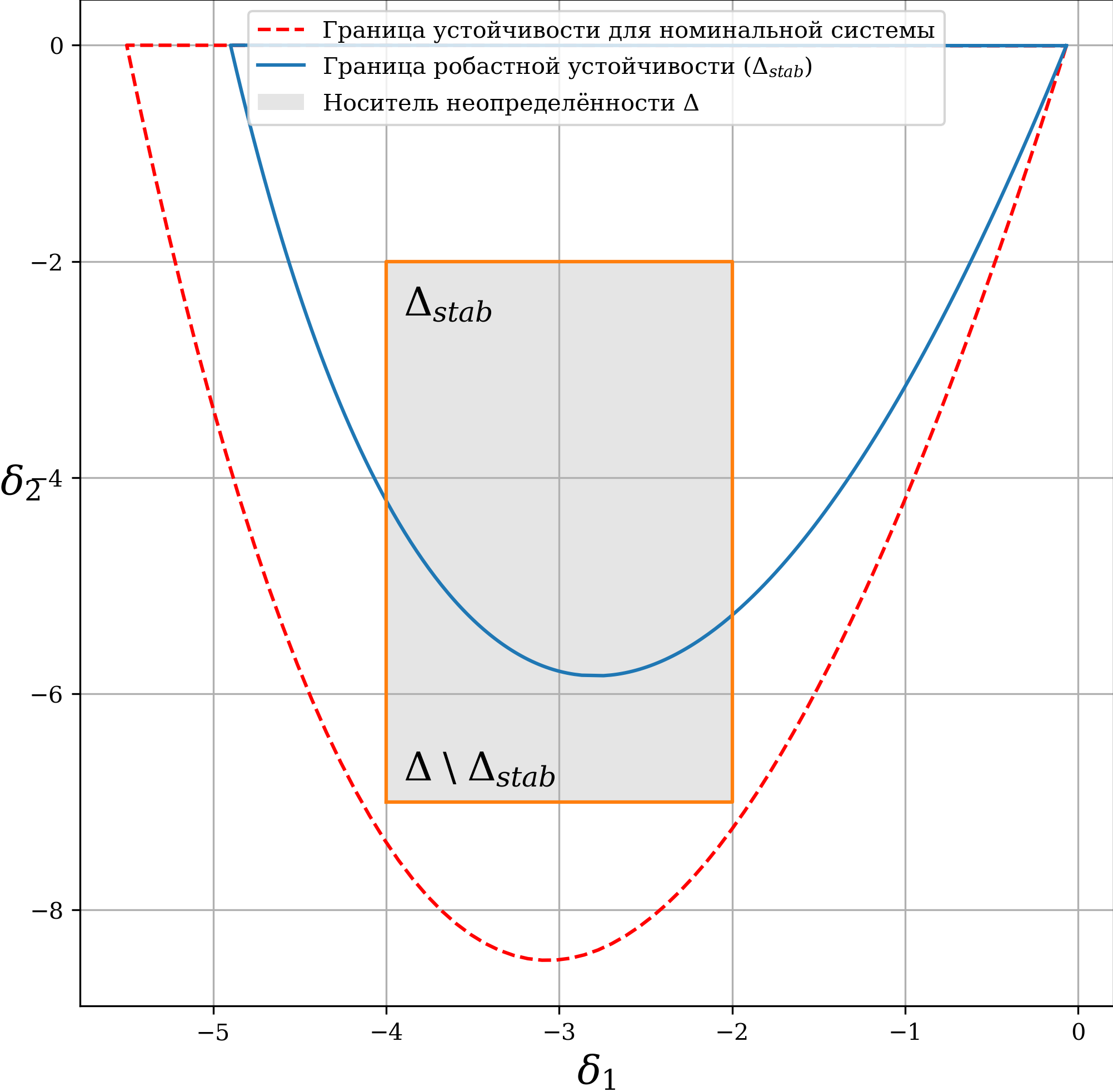}}
\caption{Номинальное и робастное $D$-разбиение для характеристического полинома
\eqref{eq:mixed-polynomial}.
Штриховая линия --- граница номинального $D$-разбиения.
Сплошная линия --- граница робастного $D$-разбиения.
Серый прямоугольник --- носитель ($\Delta$) равномерно распределенной неопределённости
$\delta$, разделённый на множество ``хороших'' ($\Delta_{stab}$) 
и ``плохих'' ($\Delta \setminus \Delta_{stab}$) параметров.}
\label{fig:2}
\end{figure}

Рассмотрим пример с замкнутой системой, имеющей характеристический
полином
\begin{equation} \label{eq:mixed-polynomial}
P(s, q, \delta) = \big((1 + q_1) s^2 + (1 + q_2) s - 1\big) s + 
(s - 15) (\delta_1 s + \delta_2 - 0{,}3 s^2).
\end{equation}
Детерминированная неопределённость $q$ 
ограничена эллипсом $100 q_1^2 + 25 q_2^2 \leq 9$,
в то время как случайный неопределённый параметр $\delta$
равномерно распределён в заштрихованном прямоугольнике
(см. рис.~\ref{fig:2}): $\delta_1 \sim [-4, -2], \delta_2 \sim [-7, -2]$.
Этот параметр представляет неопределённость коэффициентов ПИ-регулятора.
Цель состоит в решении задачи \eqref{eq:mixed-robustness-independent} 
о смешанной робастности при данных независимых неопределённостях.

Штриховая линия на рис.~\ref{fig:2} даёт границу множества устойчивости
для номинальной системы ($q_1 = q_2 = 0$). 
Поскольку все возможные случайные параметры попадают в это множество
устойчивости, то номинальная система будет робастно устойчива 
с вероятностью единица по отношению к неопределённости $\delta$.

При наличии эллиптически ограниченной неопределённости $q$
область устойчивости $\Delta_{stab}$ уменьшается, так как она
соответствует \emph{робастной устойчивости}.
Следуя первому шагу алгоритма~\ref{alg:independent-Q-Delta} из раздела~\ref{sec:main-results}, 
это множество построено с помощью робастного $D$-разбиения. 
Границы $\Delta_{stab}$ изображены на рис.~\ref{fig:2} сплошной линией.
Видно, что существуют значения неопределённости 
$\delta \in \Delta \setminus \Delta_{stab}$ вне 
этой области.

Для равномерного распределения второй шаг алгоритма~\ref{alg:independent-Q-Delta}
становится задачей о вычислении геометрической вероятности.
В итоге решение задачи о смешанной робастности 
сводится к вычислению площадей \ref{eq:probability-by-ratio}.
В этом примере 
вероятность смешанной робастности равна $p^* = 0{,}68901$. 

Для сравнения эта вероятность была оценена с помощью сценарного подхода:
для точности $\varepsilon = 0{,}01$ и уровня достоверности $\theta = 10^{-7}$ 
на множестве $\Delta$ согласно равномерному распределению и оценке \eqref{eq:chernoff-formulae} 
были выбраны $84507$ образцов случайных параметров $\delta^i$.
Для каждого образца была проверена робастная устойчивость полинома $P(s, q, \delta^i)$.
В 57856 случаях полином был робастно устойчив, что приводит к оценке
$\widehat{p}^* = 0{,}68830$, лежащей в пределах точности $|p^* - \widehat{p}^*| = 7{,}2\cdot10^{-4} < \varepsilon$.

\subsection{Скалярные случайные и детерминированные неопределённости}  \label{sec:example-3}
Рассмотрим линейную систему с характеристическим полиномом, зависящую от
скалярных неопределённых параметров $q$ и $\delta$:
\begin{equation} \label{eq:example3}
P(s, q, \delta) = s^3 + (2 - |q - \delta|) s^2 + 2 s^2 + 3,
\end{equation}
Полином устойчив при $|q - \delta| < 0{,}5$, что геометрически соответствует 
полосе между прямыми $q - \delta = 0{,}5$ и $q - \delta = -0{,}5$, 
см. рис.~\ref{fig:ex3-definition-and-Q-Delta} слева.

Далее приведены различные связи между этими двумя параметрами,
демонстрирующие все варианты задач о смешанной робастности 
\eqref{eq:mixed-robustness-independent}, 
\eqref{eq:mixed-robustness-Q(Delta)}, 
\eqref{eq:mixed-robustness-Delta(Q)}.

\subsubsection{Независимая смешанная робастность ($Q$--$\Delta$ робастность)}
\label{sec:Q-Delta-robustness}

\begin{figure} 
\centerline{\includegraphics[height=5cm]{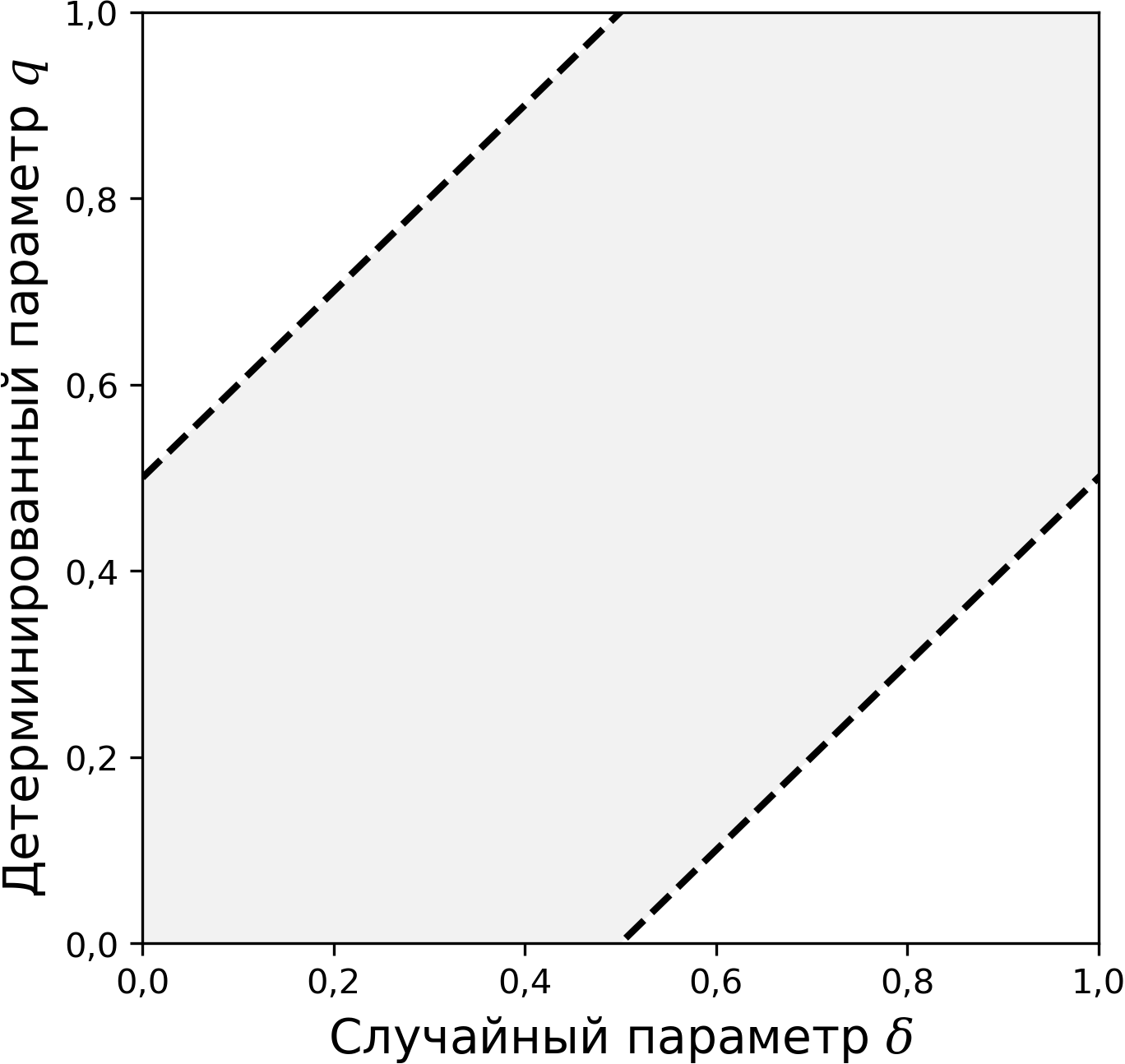}~\
\includegraphics[height=5cm]{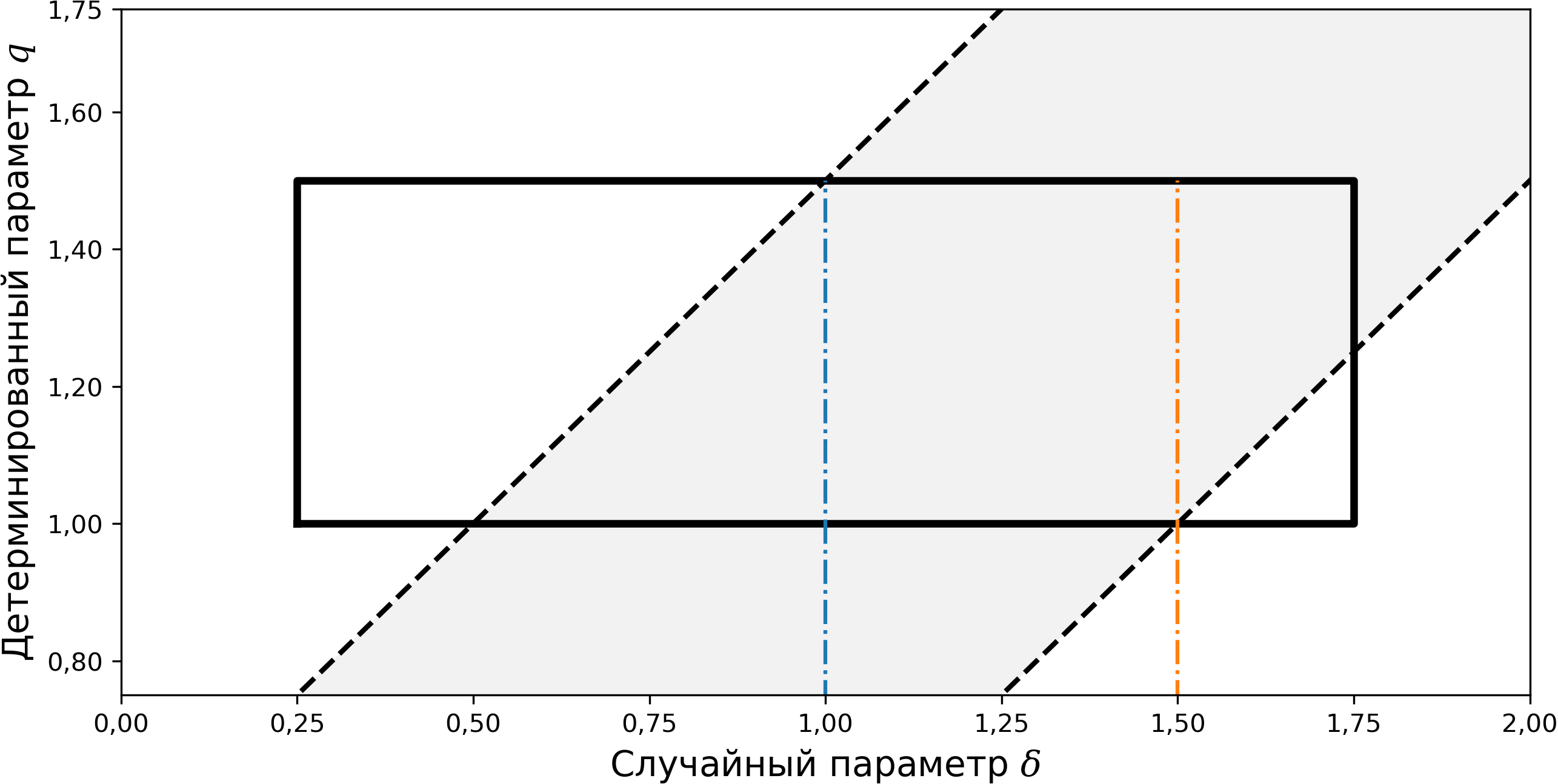}}
\caption{Пример 3. Серые области с пунктирным границами --- множества устойчивости.
Слева: область устойчивости полинома \eqref{eq:example3}.
Справа: иллюстрация решения задачи из подраздела~\ref{sec:Q-Delta-robustness} 
о $Q$--$\Delta$ робастности с равномерным распределением случайного параметра.
Штрихпунктирные линии по пересечению границ множества неопределённых параметров 
и множества устойчивости определяют подходящий отрезок $[1, 1{,}5]$ для случайного
параметра.}
\label{fig:ex3-definition-and-Q-Delta}
\end{figure}

Рассмотрим задачу о ($Q$--$\Delta$ робастности \eqref{eq:mixed-robustness-independent}.
Пусть множество значений детерминированного параметра задано интервалом $Q = [1,\, 1{,}5]$, 
а случайный параметр равномерно распределён на отрезке 
$\Delta = [0{,}25, \, 1{,}75]$, см. рис.~\ref{fig:ex3-definition-and-Q-Delta}, справа.

Здесь применим двухэтапный алгоритм~\ref{alg:independent-Q-Delta} 
из раздела~\ref{sec:main-results}.
Полином \eqref{eq:example3} робастно устойчив (по $q \in Q$) для всех 
$\delta \in \Delta_{stab} = [1, 1{,}5]$. Вне этого отрезка существуют
значения $q \in Q$, нарушающие устойчивость.
Согласно \eqref{eq:probability-by-good-set}, \eqref{eq:probability-by-ratio} 
вероятность смешанной  $Q$--$\Delta$ робастности вычисляется в явном виде
как отношение длин отрезков $\Delta_{stab}$ и $\Delta$:
$$
p^* = \frac{1{,}5 - 1}{1{,}75 - 0{,}25} = \frac{1}{3}.
$$

\subsubsection{Смешанная робастность с зависимостью первого типа ($Q(\Delta)$ робастность), первый вариант}
\label{sec:Q(Delta)-robustness-case1}

\begin{figure} 
\centerline{\includegraphics[height=6cm]{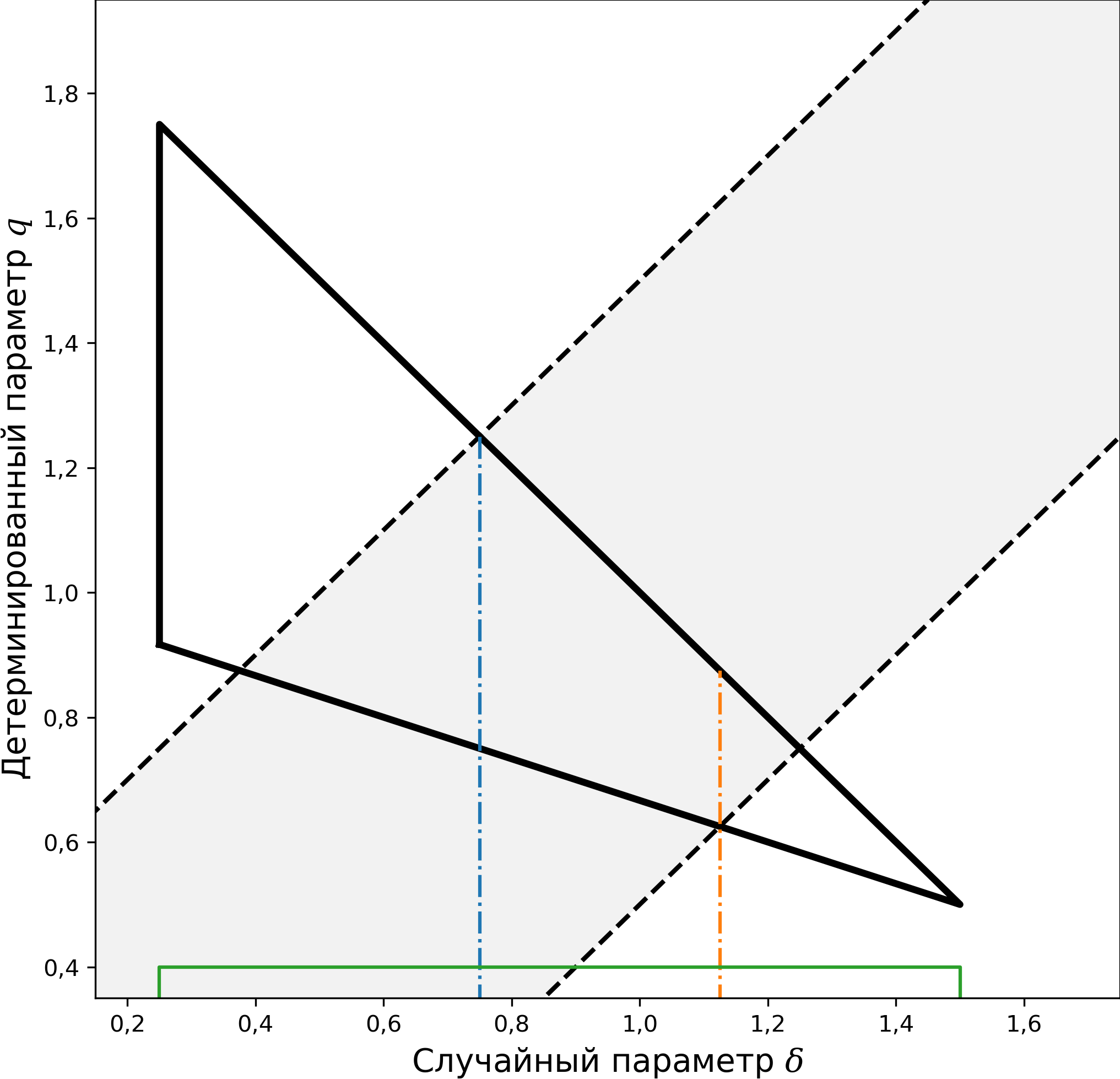}~\
\includegraphics[height=6cm]{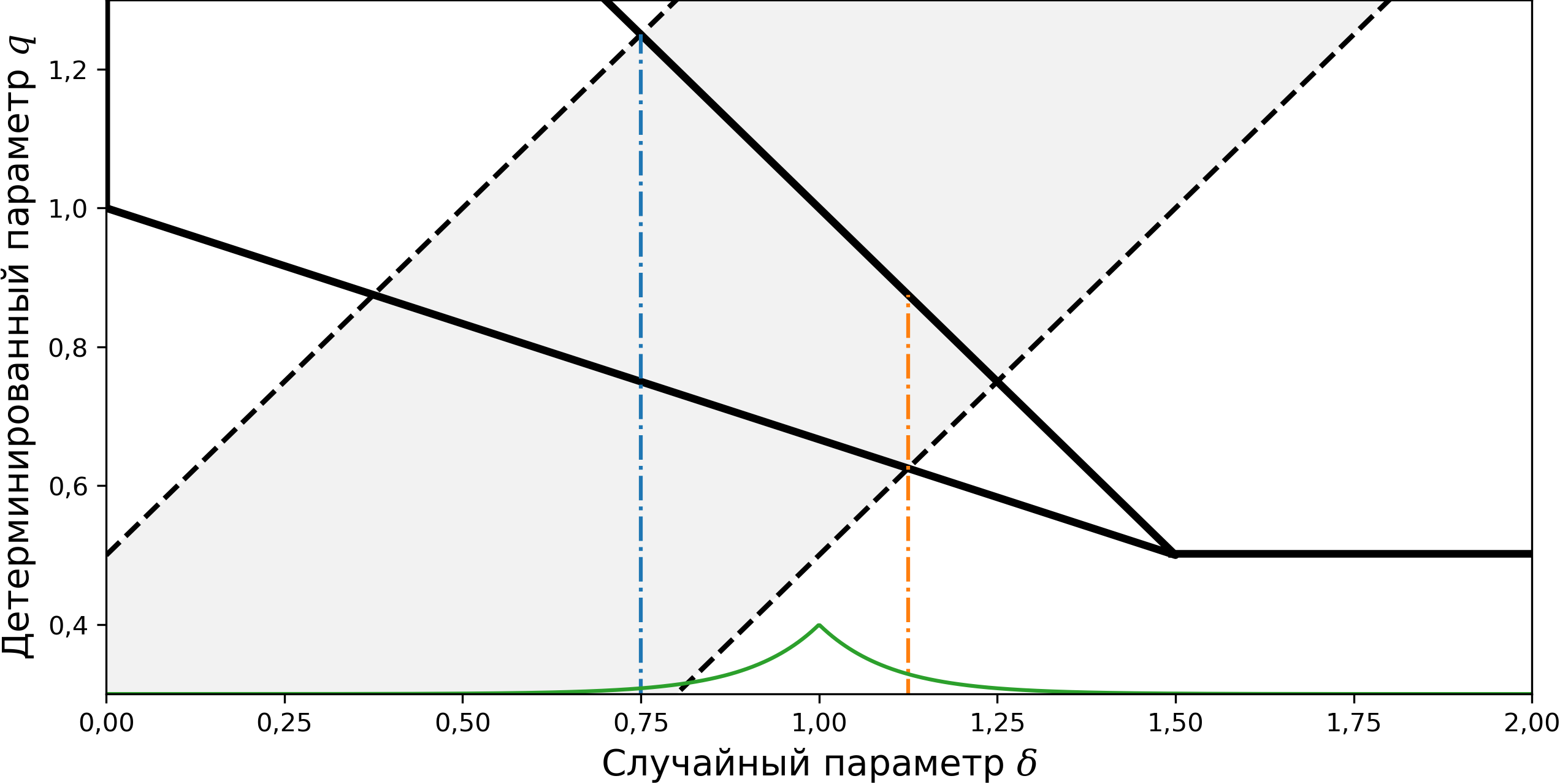}}
\caption{Серые области с пунктирным границами --- множества устойчивости.
Сплошные линии --- границы множества неопределённости.
Штрихпунктирные линии определяют края отрезков случайных параметров 
(границы множеств $\Delta_{stab} = [0{,}75, 1{,}125]$, одинаковых в обоих случаях).
Схематически обозначены плотности распределения случайного параметра.
Слева: пример~\ref{sec:Q(Delta)-robustness-case1}, 
первый вариант $Q(\Delta)$ робастности с равномерным распределением
случайного параметра.
Справа: пример~\ref{sec:Q(Delta)-robustness-case2}, второй вариант $Q(\Delta)$ робастности
с распределением Лапласа случайного параметра.}
\label{fig:ex3-Q(Delta)}
\end{figure}

Пусть случайный параметр $\delta$ равномерно распределён на отрезке $\Delta = [0{,}25, 1{,}5]$.
Множество значений детерминированного параметра представляет
собой случайное множество $Q(\delta) = [1 - \delta/3, 2 - \delta]$, 
см. рис.~\ref{fig:ex3-Q(Delta)}, слева.

В данном случае можно также применить модифицированный алгоритм~\ref{alg:independent-Q-Delta}, 
описанный в конце раздела~\ref{sec:main-results}.
Полином \eqref{eq:example3} робастно устойчив (по $q \in Q(\delta)$) 
при значениях $\delta \in \Delta_{stab} = [0{,}75, 1{,}125]$.
Как и в примере из подраздела~\ref{sec:Q-Delta-robustness}, 
вероятность смешанной $Q(\Delta)$ робастности
вычисляется как отношение отрезков
$$
p^* = \frac{1{,}125 - 0{,}75}{1{,}5 - 0{,}25} = 0{,}3.
$$

\subsubsection{Смешанная робастность с зависимостью первого типа ($Q(\Delta)$ робастность), второй вариант}
\label{sec:Q(Delta)-robustness-case2}

В этом варианте случайный параметр $\delta$ имеет распределение
Лапласа с коэффициентом сдвига $1$ и коэффициентом масштаба $0{,}1$. 
Носитель $\Delta = \mathbb{R}$ не ограничен. 
Множество значений детерминированного параметра представляет
собой случайное множество: отрезок $Q(\delta) = [1 - \delta/3, 2 - \delta]$,
если $\delta < 1{,}5$, или одну точку $Q(\delta) = \{0{,}5\}$ для $\delta \geq 1{,}5$,
см. рис.~\ref{fig:ex3-Q(Delta)}, справа.

Двухэтапный алгоритм~\ref{alg:independent-Q-Delta} применим и в этом случае, 
множество $\Delta_{stab}$ то же
и по построению не зависит от распределения.
Для скалярной случайной неопределённости и известного закона распределения
 можно использовать \eqref{eq:exact-solution-1d}:
$$
p^* 
= \mathbf{CDF}_{Laplace}(1{,}125) - \mathbf{CDF}_{Laplace}(0{,}75) 
= 1 - \frac{1}{2}e^{-\frac{1{,}125 - 1}{0{,}1}} - \frac{1}{2}e^{-\frac{1 - 0{,}75}{0{,}1}}
= 0{,}8157.
$$

\subsubsection{Смешанная робастность с зависимостью второго типа ($\Delta(Q)$ робастность)}
\label{sec:Delta(Q)-robustness}

\begin{figure} 
\centerline{\includegraphics[height=8cm]{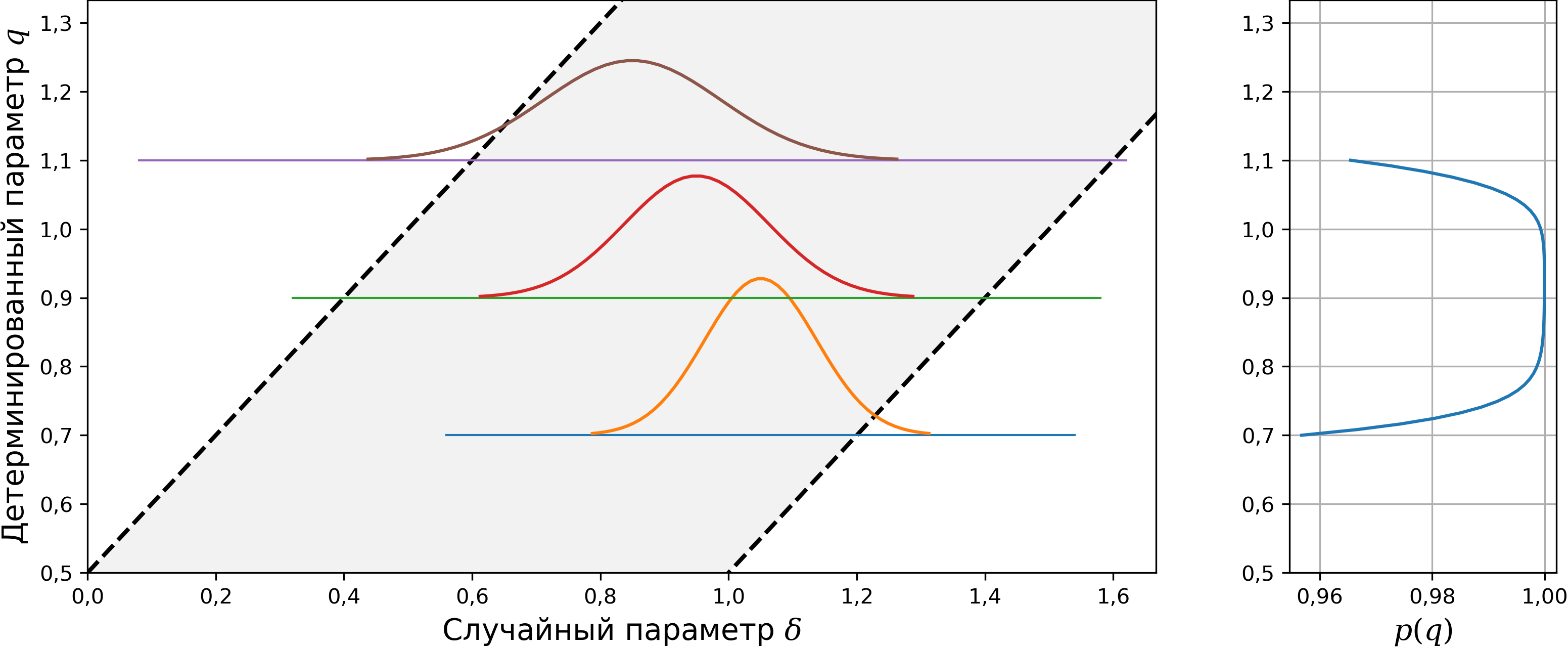}}
\caption{Пример \ref{sec:Delta(Q)-robustness}. 
Слева: пример $\Delta(Q)$ робастности 
с нормальным распределением случайного параметра.
Серая область --- множество устойчивости.
Схематически проиллюстрированы три плотности распределения
случайной величины $\delta$ для разных значений параметра $q = 1{,}1,\, 0{,}9,\, 0{,}7$.
Справа: связь между вероятностью устойчивости $p(q)$ (ось абсцисс) и параметром $q$ (ось ординат).}
\label{fig:ex3-Delta(Q)}
\end{figure}

Рассмотрим детерминированный параметр $q \in Q = [0{,}7, 1{,}1]$
и случайный параметр $\delta$ c нормальным распределением $\mathcal{N}(\delta_0, \sigma)$.
Параметры распределения зависят от детерминированного параметра:
среднее $\delta_0 = 1{,}4 - 0{,}5 q$, среднеквадратичное отклонение $\sigma = q/8$, 
см. рис.~\ref{fig:ex3-Delta(Q)}.
Отметим, что параметр $q$ не только влияет на распределение, но и входит 
в характеристический полином \eqref{eq:example3} явно.

Вероятность устойчивости $\mathrm{Prob}_{\delta \sim \Delta(q)}\,[\text{Полином } 
\eqref{eq:example3} \text{ гурвицев}]$ можно найти в явном виде по формуле
\eqref{eq:exact-solution-1d}, учитывая, что области устойчивости
задаются отрезками $D_{stab}(q) = [q - 0{,}5, q + 0{,}5]$:
\[
p(q) = \mathrm{Prob}_{\delta \sim \Delta(q)} = 
\mathop{\int}_{q - 0{,}5}^{q + 0{,}5} e^{\frac{\delta - (1{,}5 - 0{,}5 q)}{2 (q / 8)^2}} d \delta.
\]
График этой вероятности приведён на рис.~\ref{fig:ex3-Delta(Q)} справа.
Решением задачи о смешанной $\Delta(Q)$ робастности является 
минимальное значение вероятности устойчивости $p = 0{,}9568$ 
при ``наихудшем'' $q = 0{,}7$.

\section{Заключение и расширения}
\label{sec:conclusion}

В статье проведён обзор различных типов параметрической робастности.
На примере исследования устойчивости систем с детерминированной и случайной
неопределённостями предложены постановки задачи о смешанной робастности.

Задачи отличаются взаимозависимостью детерминированных и случайных
групп параметров. Заменой критерия эти постановки элементарно обобщаются на другие свойства
системы, например на (робастную) апериодичность и пр.
Можно исследовать иные зависимости детерминированных и случайных 
параметров, в том числе их иерархические графы.

Предлагается алгоритм, разделяющий решение основной задачи
о смешанной робастности на два независимых этапа.
На первом этапе метод анализа в пространстве параметров 
описывает подмножество ``хороших'' случайных параметров.
На втором этапе используется другой метод, который аналитически, численно или методом 
Монте-Карло оценивает вероятность попадания
случайного параметра в указанное множество.

Используемый на первом этапе метод робастного $D$-разбиения
позволяет решать и другие задачи о смешанной робастности.

Во-первых, помимо устойчивости
можно проверять 
нечувствительность к параметрическим возмущениям заданных показателей качества,
например запаса устойчивости по усилению и фазе, 
коэффициента демпфирования и т.п. Эти возможности открываются
благодаря тому, что подобные показатели 
определяются расположением
корней характеристического полинома или характеризуются
с помощью частотных методов.
Одновременная робастность по отношению к нескольким критериям качества 
исследуется аналогично.

Во-вторых, наряду с линейными системами непрерывного времени 
предложенный подход позволяет исследовать смешанную робастную устойчивость 
иных типов линейных систем с постоянными параметрами.
Так, устойчивость дискретных систем определяется свойством Шура характеристического
полинома;
устойчивость систем с запаздыванием определяется свойством Гурвица
характеристического квазиполинома;
достаточное условие устойчивости систем дробного порядка
также описывается расположением корней 
``полинома'' с нецелыми степенями.

Кроме точных решений, в подразделе~\ref{sec:approximations}
предложен ряд процедур, позволяющих решать
задачи о смешанной робастности приближённо.

Помимо решённых, поставлен ряд открытых вопросов:
как эффективно решать задачи \eqref{eq:mixed-robustness-Q(Delta)} 
и \eqref{eq:mixed-robustness-Delta(Q)} 
о смешанной робастности 
с зависимыми детерминированными и случайными параметрами?
Насколько хорошо их удаётся решать приближенными методами?
На примере из подраздела~\ref{sec:example-3}) показано, что 
решение задач смешанной робастности упрощается, если
известно подмножество $Q \times \Delta$ объединённого набора
параметров, на котором система устойчива. 
Как описать это множество для иных свойств системы, аналогично тому, 
как аппарат $D$-разбиения используется для устойчивости?

Автор выражает сердечную признательность Б.Т.~Поляку за поддержку при подготовке статьи,
а также благодарит анонимного рецензента, замечания которого позволили существенно улучшить 
статью.

\end{document}